\documentclass[lettersize,journal]{IEEEtran}
\usepackage{amsmath,amsfonts}
\usepackage{algorithmic}
\usepackage{array}
\usepackage[caption=false,font=normalsize,labelfont=sf,textfont=sf]{subfig}
\usepackage{textcomp}
\usepackage{stfloats}
\usepackage{url}
\usepackage{verbatim}
\usepackage{graphicx}
\usepackage{cite}
\newtheorem{lemma}{Lemma}[section]
\newtheorem{prop}{Proposition}[section]
\newtheorem{te}{Theorem}[section]

\newtheorem{ass}{A}

\newcommand{\be}{\begin{equation}}
\newcommand{\ee}{\end{equation}}
\newcommand{\ba}{\begin{array}}
\newcommand{\ea}{\end{array}}
\newcommand{\bee}{\begin{eqnarray*}}
\newcommand{\eee}{\end{eqnarray*}}
\newcommand{\bea}{\begin{eqnarray}}
\newcommand{\eea}{\end{eqnarray}}

\newcommand{\II}{\mathbb{I}}
\newcommand{\WW}{\mathbb{W}}
\newcommand{\TT}{\mathbb{T}}
\newcommand{\HH}{\mathbb{H}}
\newcommand{\RR}{\mathbb{R}}
\newcommand{\BB}{\mathbb{B}}
\newcommand{\CC}{\mathbb{C}}

\newcommand{\DD}{\mathbb{D}}

\newcommand{\GG}{\mathbb{G}}

\newcommand{\mx}{\mathbf{x}}
\newcommand{\mz}{\mathbf{z}}

\newcommand{\mv}{\mathbf{v}}
\newcommand{\mq}{\mathbf{q}}
\newcommand{\mg}{\mathbf{g}}
\newcommand{\mc}{\mathbf{c}}
\newcommand{\uu}{\mathbf{u}}
\newcommand{\md}{\mathbf{d}}
\newcommand{\mr}{\mathbf{r}}
\newcommand{\me}{\mathbf{e}}

\begin{document}

\title{A Hessian inversion-free exact second order method for distributed consensus optimization}

\author{Du\v{s}an Jakoveti\'c,~\IEEEmembership{Member,~IEEE,} Nata\v sa Kreji\'c, Nata\v sa Krklec Jerinki\'c  
\thanks{The authors are with the University of Novi Sad, Faculty of Sciences, Department of Mathematics and Informatics, Trg Dositeja Obradovica 4, Novi Sad, Serbia. Authors' emails: dusan.jakovetic@dmi.uns.ac.rs; natasak@uns.ac.rs; natasa.krklec@dmi.uns.ac.rs.    
This work is supported by the  Ministry of Education,
Science and Technological Development, Republic of Serbia. It is also supported in part by the Bilateral Cooperation Serbia – Croatia, project ``Optimization Method Application in Biomedicine,'' and by the European Union’s Horizon 2020 Research and Innovation program under grant agreement No 871518. The paper reflects only the view of the authors and the Commission is not responsible for any use that may be made of the information it contains.}
}

\markboth{Journal of \LaTeX\ Class Files,~Vol.~14, No.~8, August~2021}%
{Shell \MakeLowercase{\textit{et al.}}: A Sample Article Using IEEEtran.cls for IEEE Journals}


\maketitle

\begin{abstract}
We consider a standard distributed consensus optimization problem where a set of agents connected over an undirected network minimize the sum of their individual (local) strongly convex costs. Alternating Direction Method of Multipliers (ADMM) and Proximal Method of Multipliers (PMM) have been proved to be effective frameworks for design of exact distributed second order methods (involving calculation of local cost Hessians). However, existing methods involve explicit calculation of local Hessian inverses at each iteration that may be very costly when the dimension of the optimization variable is large. In this paper, we develop a novel method, termed INDO (Inexact Newton method for Distributed Optimization), that alleviates the need for Hessian inverse calculation. INDO follows the PMM framework but, unlike existing work, approximates the Newton direction through a generic fixed point method (e.g., Jacobi Overrelaxation) that does not involve Hessian inverses. We prove exact global linear convergence of INDO and provide analytical studies on how the degree of inexactness in the Newton direction calculation affects the overall method’s convergence factor. Numerical experiments on several real data sets demonstrate that INDO’s speed is on par (or better) as state of the art methods iteration-wise, hence having a comparable communication cost. At the same time, for sufficiently large optimization problem dimensions $n$ (even at $n$ on the order of couple of hundreds), INDO achieves savings in computational cost by at least an order of magnitude.
\end{abstract}

\begin{IEEEkeywords}
Inexact Newton, proximal method of multipliers, distributed optimization, exact convergence, strongly convex problems.
\end{IEEEkeywords}

\section{Introduction}

We consider problems of the form
\be \label{51} \min_{y \in \mathbb{R}^n} \sum_{i=1}^{N} f_i(y).
\ee
Here, $f_i : \mathbb{R}^n \to \mathbb{R}$, $i=1,...,N$, is a strongly convex local cost function assigned to a node within a network of distributed agents able to perform local operations and communicate with their neighbours. 
 Formulation~(1) finds a number of applications in signal processing, e.g.,  \cite{novo1,SayedEstimation}, control, e.g.,  \cite{JoaoMotaMPC}, Big Data analytics, e.g.,  \cite{scutari2}, social networks, e.g.,  \cite{novo3}, etc. The available methods for solving (\ref{51}) include a large class of the so called exact methods that ensure convergence to the exact solution of (\ref{51})  with different rates of convergence.  The exact convergence is achieved in several ways -- by utilizing diminishing step sizes in gradient methods for penalized reformulation \cite{arxivVersion,Nedickaskadni}, by gradient tracking or second order  methods that are defined within primal-dual framework, e.g,  \cite{Espectral, harnessing, dusan,ESOM,extra,EISEN,UsmanXin,d3,d5,d6,optra,d8,dragana,Zhang,Sayed2,Sayed3},  or in the framework of alternating direction methods \cite{BoydADMM, DQM}. 
  Multiple consensus steps per each gradient update to ensure exact convergence are also  considered~\cite{wei}.

For the current paper, of special relevance are two strategies available in the literature -- the proximal method of multipliers as a framework to develop exact distributed methods in \cite{ESOM}, and the well known theory of inexact Newton methods in centralized optimization, \cite{DES}.  The main advantage of the inexact Newton methods is that they avoid oversolving of the Newtonian system of linear equations; as such, they have been employed in related problems, for example in minimizing finite sums in machine learning applications  \cite{IMA}. To be more specific, we consider the constrained reformulation of (\ref{51}) in the augmented space (like in, e.g., \cite{ESOM}) and build up on the distributed second  order approximation of the Augmented Lagrangian. The second order approximation of the Augmented Lagrangian conforms with the sparsity structure of the network and hence Newton-like step is well defined. The main obstacle for applying the Newton method is the fact that although the Hessian is sparse and distributed, the inverse Hessian is dense and challenging to compute efficiently in a distributed environment. One possibility is to approximate the inverse Hessian by inverting (local) Hessians of the $f_i$'s and build a Taylor-like  approximation to the inverse of the global Hessian as suggested in \cite{NN}, \cite{ESOM}, and exploited in \cite{r1,r2}.

The approach we propose here is based on a recent result on distributed solution of systems of linear equations \cite{DFIX}. More precisely, the Newton-like equation that defines the update step is a solution of a system of linear equations defined by the global Hessian of the Augmented Lagrangian. 
So, instead of approximating the inverse Hessian as done in, e.g., \cite{NN,ESOM}, we propose to solve  the system of linear equations inexactly, in a distributed manner by a suitable iterative method. 
 A distinctive feature of the approach is that it avoids inversion of the Hessians of the $f_i$'s; instead, only diagonal elements of the Hessians are inverted. Another appealing feature of the approach is compatibility with  the theory of inexact Newton methods in centralized optimization. The proposed approach  is particularly suitable for the case of relatively large $ n, $ when the approximation of the inverse Hessian as in \cite{NN, ESOM} might be expensive, and thus an iterative solver for the system of linear equations is a natural choice.  The convergence theory developed here relies   on the theory for the proximal methods of multipliers, as the proposed method fits this general framework, following, e.g., \cite{ESOM}. However, the error analysis carried out here -- reflecting the inexact solution of the Newtonian system of linear equations and how it affects the overall proximal multipliers convergence -- is very different for the proposed method.  

To summarize, the main contributions of this paper are the definition and convergence theory of the novel exact second order method termed INDO (Inexact Newton method for Distributed Optimization).\footnote{To avoid confusion, we clarify that the wording ``inexact Newton'' here refers to approximately solving the Newtonian system of linear equations that defines the Newton direction; on the other hand, the wording ``exact method'' here signifies that INDO converges to the exact solution of~(1).}  Linear convergence to the exact solution is shown for strongly convex problems under a set of standard assumptions. The rate of convergence is also analysed, assuming that the linear solver is of fixed-point type, in particular we focused on the application of the Jacobi Overrelaxation method to derive an estimate of the convergence factor. The convergence factor is compared with the factor corresponding to ESOM method, \cite{ESOM}, as ESOM operates within the same framework -- proximal method of multipliers with approximate second order direction. Furthermore, we analyse computational and communication costs of these two methods and show that the computational costs of INDO can be significantly smaller than the costs of ESOM while the communication costs per iteration are the same. These claims are illustrated by numerical results for  quadratic and logistic loss problems including those from standard machine learning real data test cases. 

This paper is organized as follows. In Section 2 we introduce the problem model, the primal-dual framework based on the proximal method of multipliers and state the assumptions on the problem. The INDO method is introduced and theoretically analysed in Section 3. Details on distributed implementation as well as the analysis of the convergence factor are presented in Section 4, while Section 5 is devoted to practical implementation of INDO with the analysis of computational costs and comparison with ESOM \cite{ESOM}  method through  standard examples from machine learning.  Some conclusions are drawn in Section 6. Finally, lengthy supporting proofs are delegated to the Appendix.

\section{Preliminaries}

The notation used is the paper is the following. We use upper case blackboard bold letters to denote matrices in $  \mathbb{R}^{nN \times nN}  $, e.g.,  $ \mathbb{A}, \mathbb{B}, \ldots $. For $\mathbb{A} \in \mathbb{R}^{nN \times nN}$, we use its block elements representation $ \mathbb{A} = [A_{ij}], \; A_{ij} \in \mathbb{R}^{n \times n} $. 
We denote scalar elements of $\mathbb{A}$ by  $ a_{ij} \in \mathbb{R}. $  Similarly,  upper normal letters $ A,B,...$, are used for stand-alone matrices in $ \mathbb{R}^{n \times n}.  $ The vectors in $\mathbb{R}^{nN}$ are denoted by bold lowercase letters, e.g., $ \mathbf{x} \in \mathbb{R}^{nN}  $; their component blocks are $  x_i \in \mathbb{R}^n $; similarly, we use normal lowercase letters, e.g., $y$, for stand-alone vectors in $ \mathbb{R}^n. $ The notation $ \|\cdot\|_2 $ stands for the 2-norm of its vector or matrix argument; if the subscript is omitted,  $\|\cdot\|$ also designates the 2-norm, if not stated otherwise.  We denote by $diag(A)$ the diagonal matrix with diagonal elements equal to those of matrix $A$, 
and by $\sigma(A)$ the spectral radius of~$A$.

The network of connected  computational nodes (agents)  is  represented by a graph $G=(V, \mathcal{E})$, where $V$ is the set of nodes $\{1,...,N\}$ and $\mathcal{E}$ is the set of undirected edges $(i,j)$.  Denote by $O_i$ the set of neighbors of node $i$ and let $\bar{O_i}=O_i \bigcup \{i\}.$  The communication network is accompanied by a communication (weight) matrix $ W $ with the following properties. 

 \begin{ass} \label{A1}  The matrix $ W  \in \mathbb{R}^{N \times N} $  is  symmetric, doubly stochastic
 and 
 $$ w_{ij} > 0 \mbox{ if } j\in \bar{O}_i, \;  w_{ij} = 0 \mbox{ if }  j \notin \bar{O}_i $$  \end{ass}

 Let us assume that each of the $N$ nodes has its local cost function $f_i$ and has access to  computing its first and second derivatives. Under the assumption A\ref{A1}, the problem \eqref{51} has the equivalent form 
\be \label{52} \min_{\mathbf{x} \in \mathbb{R}^{nN}} f(\mathbf{x}):=\sum_{i=1}^{N} f_i(x_i) \quad \mbox{s. t. } \quad (\II - \WW)^{1/2} \mathbf{x}=0,
\ee
where $\mathbf{x}=(x_1;...;x_N) \in \RR^{nN}$, $\WW=W \otimes I \in \RR^{nN \times nN}$ and $\II \in \RR^{nN \times nN} $ is the identity matrix. 

\begin{ass} \label{A2} The  functions $ f_i $'s are twice continuously differentiable and the eigenvalues of the local Hessians are bounded by positive constants $ 0 < m \leq  M < \infty, $ i.e.,
 \be
\nonumber
m I \preceq \nabla ^2f_i(y) \preceq M I,
\ee
for all $ y \in \mathbb{R}^n $ and $ i=1,\ldots,N. $ 
\end{ass}

The above assumption implies that the   function $ f(\mx) $  is also strongly convex  with constant $ m $ and its gradient $ \nabla f$ is Lispchitz continuous with constant $ M. $ Given that we consider a method of Newton-type, the assumption on Lipschitz continuous Hessian is also needed, as usual in the theory of Newton methods and in \cite{ESOM}. 

\begin{ass} \label{A3} The Hessian of the objective function $ \nabla^2 f(\mx) $ is Lipschitz continuous, i.e., there exists $ L > 0 $ such that
\be 
\nonumber
\|\nabla^2f(\mx) - \nabla^2f(\mz)\| \leq L \|\mx-\mz\|, \; \; \mx,\mz \in \mathbb{R}^{nN}. 
\ee
\end{ass}

The following Lemma from \cite{DQM} will be used in the convergence analysis. 

\begin{lemma} \cite{DQM}  \label{lemaM} Assume that A\ref{A2}-A\ref{A3} hold. Then for all $ \mx,\mz \in \mathbb{R}^{nN} $ 
\begin{eqnarray*}
 &\,& \|\nabla f(\mx) - \nabla f(\mz) + \nabla^2 f(\mx)(\mz-\mx)\| \\
 &\,&\leq  \min\{2M,\frac{L}{2}\|\mx-\mz\|\}\|\mx-\mz\|. 
 \end{eqnarray*}
\end{lemma}

Now, 
the Augmented Lagrangian function of (\ref{52}) is defined as 
\be \nonumber
{\cal L}(\mx,\mv) = f(\mx) + {\mv}^T(\II-\WW)^{1/2}\mx + \frac{\alpha}{2} \mx^T(\II-\WW)\mx,
\ee where $ \alpha > 0 $ is a constant and $ {\mv} $ is  dual variable. For a strictly positive proximal coefficient $ \varepsilon, $ the proximal method of multipliers is defined by the primal update
\be
\label{primal} \mx^{k+1} = \arg \min {\cal L}(\mx,{\mv}^k) + \frac{\varepsilon}{2}\|\mx-\mx^k\|^2, \ee
while the dual update with the stepsize $ \alpha $  is 
\be \label{dual} {\mv}^{k+1} = {\mv}^k + \alpha (\II -\WW)^{1/2} \mx^{k+1}. \ee
The primal step is not implementable in the distributed environment due to the augmented term $ \frac{\alpha}{2} \mx^T(\II-\WW)\mx $ and therefore an approximation of $ {\cal L}(\mx,\mv) $ is needed. We consider the second order Taylor approximation with respect to $ \mx $ at the point $ (\mx^k, {\mv}^k), $
\be 
\label{54} 
{\cal L}(\mx,{\mv}^k) \approx {\cal L}(\mx^k,{\mv}^k) + \mg^k (\mx-\mx^k) + \frac{1}{2} (\mx-\mx^k)^T \HH^k (\mx-\mx^k), 
\ee
with \be
\label{55}
\HH^k = \nabla^2 f(\mx^k) + \alpha(\II-\WW) + \varepsilon \II. 
\ee
and 
\be
\nonumber
{\mg}^k:=\nabla_{\mx} {\cal L}(\mx^k,{\mv}^k)  = \nabla f(\mx^k) + (\II-\WW)^{1/2} {\mv}^k + \alpha(\II - \WW)\mx^k. 
\ee
Using the second order Taylor approximation (\ref{54}) one can define the Newton step $ \md_{N}^k $ as
\be 
\label{57}
\HH^k \md_{N}^k = -\mg^k ,
\ee
and form the primal update as $ \mx^{k+1} = \mx^k + \md_N^k. $ Notice that the matrix $ \HH^k $ has the sparsity structure of the graph; hence, one can apply an iterative method of the fixed point-type to solve (\ref{57}). The details of such procedure can be found in \cite{DFIX}. However,  solving (\ref{57}) exactly might be too costly and an approximate solution to (\ref{57}) is in fact sufficient for the convergence as we will demonstrate further on. Thus the step we use is defined by 
$$
\HH^k \md^k = - \mg^k + \mr^k ,$$
for some residual $ \mr^k $ which obeys the classical Inexact Newton forcing condition \cite{DES}  
$$ \|\mr^k \| \leq \eta_k \|\mg^k\|, $$
for a forcing term $ \eta_k \geq 0. $ 

The dual update is given by (\ref{dual}) but the matrix $  \alpha (\II -\WW)^{1/2} $ is not neighbor sparse. Thus we apply the same variable transformation as in \cite{ESOM}, defining a sequence of variables $ \mq^k $ as 
 $$ \mq^k  =  (\II -\WW)^{1/2} {\mv}^k. $$
  Now, multiplying the dual update (\ref{dual}) by $   (\II -\WW)^{1/2} $ from the left and using the definition of $ \mq^k $ we obtain
\be 
\label{dualq}
\mq^{k+1} = \mq^k +  \alpha (\II -\WW) \mx^{k+1}, 
\ee
which is computable in the distributed manner. Notice that with this change of variables we get
\be
\label{59}
\mg^k = \nabla f(\mx^k) + \mq^k +  \alpha (\II - \WW) \mx^{k}. 
\ee

\section{Distributed Inexact Newton Algorithm}

With the notation introduced in the previous section we are now in a position to state the general INDO method. Further details regarding the inexact Newton step computation  (step S1 in the INDO method below), are postponed to Section 4. Therein, we discuss different possibilities for fulfilling the Inexact Newton condition (\ref{58}) in a distributed environment. Assume that the forcing sequence $ \{\eta_k\}_{k=0}^{\infty} $ is nonnegative.  

\noindent{\bf INDO method}
\begin{itemize}
\item[Input]: $ \mx^0 \in \mathbb{R}^{n N}, \; \mq^0 = 0, \; k=0, \{\eta_k\}_k, \alpha > 0, \; \varepsilon > 0 $
\item[S1] Find a step $ \md^k $ such that 
\be \label{58}  \HH^k \md^k =- \mg^k + \mr^k, \; \|\mr^k\| \leq \eta_k \|\mg^k\|. \ee
\item[S2] 
Compute the primal update
$$ \mx^{k+1} = \mx^k + \md^k. $$
\item[S3] Compute the dual update
$$ \mq^{k+1} = \mq^k + \alpha(\II-\WW)\mx^{k+1},$$ 
 set  $ k = k+1  $ and return to S1.
\end{itemize}

As already stated, one can easily see that the algorithm fits the general proximal multipliers framework \cite{ESOM}, with the primal step defined in S1-S2 and the dual update in S3. The key novelty is the direction computation in S1 where we generate a suitable inexact direction such that the residual is small enough with respect to the gradient norm, (\ref{58}). This important property actually allows for flexible approach in the linear solver and avoids oversolving. That is, the proposed approach allows using a relatively large tolerance while the gradient is large, and implies stringent tolerance condition as the value of gradient decreases. The details of S1 implementation are presented in Section 4, where we also give INDO algorithm from the perspective of each node in the network. 

 Due to the convexity and the fact that we only have equality constraints, $\mx^* $ is a solution of problem (\ref{52})   if and only if there exists ${\mv}^*$ such that 
\begin{eqnarray}
\label{61}
\nabla f(\mx^*) + (\II-\WW)^{1/2} {\mv}^* = 0  \\
(\II - \WW)^{1/2} \mx^* = 0. \label{62}
\end{eqnarray}
Strong convexity implies that $ (\mx^*,{\mv}^*)$  is unique. 

The convergence analysis we develop below relies on the reasoning from \cite{ESOM} in general but the error analysis, with respect to the proximal method of multipliers (see Section 2) is fundamentally different. We present the proof  through a sequence of technical Lemmas, some of which correspond to  parts of proofs in \cite{ESOM} that are rather similar. Proofs of the Lemmas similar to \cite{ESOM} are provided in Appendix for completeness. 

\begin{lemma} \label{L1} Let  A\ref{A1}-A\ref{A2}   hold. A  sequence generated by Algorithm INDO satisfies 
$$ {\mv}^{k+1} - {\mv}^{k} - \alpha(\II-\WW)^{1/2}(\mx^{k+1}-\mx^*) = 0. $$
\end{lemma}
{\em Proof. } Given the rule for dual update in S3 and the definition of $ \mq^k $ in (\ref{dualq}), we have 
$$ {\mv}^{k+1}-{\mv}^k =  \alpha(\II-\WW)^{1/2} \mx^{k+1}. $$ Adding the zero term from (\ref{62}) multiplied by $\alpha $, i.e.,   $\alpha(\II-\WW)^{1/2} \mx^* $,  implies the statement.  $ \Box $

\begin{lemma} \label{L2} Let assumptions   A\ref{A1}-A\ref{A2}  hold and $ \{\mx^k,\mq^k\} $ be a sequence generated by Algorithm INDO. Then
$$ \nabla f(\mx^{k+1})- \nabla f(\mx^*) + \varepsilon \md^k +  (\II-\WW)^{1/2}({\mv}^{k+1}-{\mv}^*) + \me^k = 0, $$
where 
$$ \me^k = \nabla^2 f(\mx^k) \md^k + \nabla f(\mx^k) - \nabla f(\mx^{k+1}) - \mr^k, \; \mr^k = \HH^k \md^k + \mg^k. $$
\end{lemma}

{\em Proof. } By the definition of step $ \md^k $ in S2, (\ref{dualq}) and (\ref{61}) we have
 
 {\allowdisplaybreaks{
 \begin{eqnarray*}
0 & = & \HH^k \md^k + \mg^k - \mr^k \\
& = & (\nabla^2 f(\mx^k) +  \alpha(\II-\WW) +\varepsilon \II)\md^k \\
&+&  \mg^k - \mr^k \pm (\II-\WW)^{1/2}{\mv}^{k+1} \\
&=& \nabla^2 f(\mx^k) \md^k + \alpha(\II-\WW)\md^k + \varepsilon \md^k + \nabla f(\mx^k) + \mq^k \\
&+& \alpha(\II-\WW)\mx^k -\mr^k \pm \mq^{k+1} \\
& = & \nabla^2 f(\mx^k) \md^k + \alpha(\II-\WW)\mx^{k+1} + \varepsilon \md^k + \nabla f(\mx^k)
\\
&+& (\II-\WW)^{1/2} {\mv}^{k+1} -\mr^k \\
&-& (\mq^{k+1}-\mq^k) - (\nabla f(\mx^*) + (\II-\WW)^{1/2} {\mv}^*) \pm \nabla f(\mx^{k+1}) \\
&=& \nabla f(\mx^{k+1}) - \nabla f(\mx^*) + \varepsilon \md^k \\
&+& (\II-\WW)({\mv}^{k+1}-{\mv}^*) + \me^k. 
\end{eqnarray*} }}
$ \Box $

\begin{lemma} \label{L3} Let assumptions A\ref{A1}-A\ref{A3} hold and $ \{\mx^k,\mq^k\} $ be the sequence generated by INDO. Then
 \begin{eqnarray*} 
 &\,& \|\me^k\| \leq \min\{2M, \frac{L}{2}\|\md^k\|\} \|\md^k\| \\
 &+& \eta_k(2\alpha+M) \|\mx^k-\mx^*\| \sqrt{2} \eta_k\|{\mv}^k-{\mv}^*\|. 
 \end{eqnarray*}
\end{lemma} 

{\em Proof. }
By Lemma \ref{lemaM}, Lemma \ref{L2} and the  condition in S2 of the algorithm we have
\begin{eqnarray} 
\|\me^k\| & \leq & \|\nabla^2 f(\mx^k) \md^k + \nabla f(\mx^k) - \nabla f(\mx^{k+1})\| + \|\mr^k\| \nonumber \\
& \leq & \min\{2M,\frac{L}{2}\|\md^k\|\}\|\md^k\| + \eta_k \|\mg^k\|. \label{ek} 
\end{eqnarray}
Furthermore, the definition of $ \mg^k $ in (\ref{59}) and optimality conditions (\ref{61}), (\ref{62}) imply
\begin{eqnarray*}
\|\mg^k\| & = & \|\nabla f(\mx^k) + (\II-\WW)^{1/2} {\mv}^k + \alpha (\II-\WW) \mx^k \\
& - & (\nabla f(\mx^*)  + (\II-\WW)^{1/2} {\mv}^* ) -  \alpha(\II-\WW) \mx^* \|\\
& \leq & \|\nabla f(\mx^k) - \nabla f(\mx^*)\| + \|(\II-\WW)^{1/2}({\mv}^k-{\mv}^*)\|\\
&+& \| \alpha(\II-\WW)(\mx^k - \mx^*)\| \\
& \leq & M\|\mx^k - \mx^*\|+ \sqrt{2} \|{\mv}^k - {\mv}^*\| \\
&+& 2 \alpha \|\mx^k-\mx^*\|\\
& = & (M + 2 \alpha) \|\mx^k - \mx^*\| + \sqrt{2} \|{\mv}^k - {\mv}^*\|. 
\end{eqnarray*}
Placing the last inequality into (\ref{ek}) we get the statement. $ \Box$ 

In the convergence analysis below we will use the following inequality, \cite{ESOM}. Let $ a,b $ be two arbitrary vectors of the same dimension and $ \xi > 0 $ be an arbitrary real number. Then
\be \label{xi}
-2a^Tb \leq \frac{1}{\xi} \|a\|^2 + \xi \|b\|^2. 
\ee
The technical lemma below is proved in \cite{ESOM}, although it is not stated separately, see the proof of Theorem 2, inequality (93) in \cite{ESOM}. 

\begin{lemma}\label{L4}
Assume that A\ref{A1}-A\ref{A3} hold and $ \{{\mx}^k,\mq^k\} $ be a sequence generated by Algorithm INDO. Then 
\begin{eqnarray*} 
&\,&\|{\mv}^{k+1}-{\mv}^*\|^2 \leq \frac{\beta \varepsilon^2}{(\beta-1)\lambda_2} \|\md^k\|^2  \\
&+& 
\frac{\phi \beta}{\lambda_2} \|\nabla f({\mx}^{k+1}) - \nabla f({\mx}^*)\|^2 + \frac{\beta \phi}{(\phi-1) \lambda_2} \|\me^k\|^2, 
\end{eqnarray*}
where $ \lambda_2 $ is the smallest nonzero eigenvalue of $ \II-\WW, $ and $ \beta, \phi > 1 $ are arbitrary constants. 
\end{lemma}

Let us define the sequence of concatenated dual and primal errors as well as the Lyapunov function using the matrix $ \cal G $ below,
$$ \uu = \begin{bmatrix} \mv \\
{\mx} \end{bmatrix} \; {\cal G} = \begin{bmatrix} \II & 0 \\ 0 & \alpha \varepsilon \II \end{bmatrix}. $$ Analogously define the concatenated $ \uu^* $ and consider the sequence $ \|\uu^k - \uu^*\|^2_{\cal G} = \|{\mv}^k - {\mv}^*\|^2 + \alpha \varepsilon \|{\mx}^k - {\mx}^*\|^2. $ We will prove that that the sequence $ \|\uu^k - \uu^*\|^2_{\cal G} $ converges to zero linearly and hence the sequence of primal errors $ \|{\mx}^{k} - {\mx}^*\| $ converges to zero   linearly  as well. 

The following Lemma is proved in the Appendix. 

\begin{lemma} \label{L5}
Assume that A\ref{A1}-A\ref{A3} hold and let $ \{{\mx}^k,\mq^k\} $ be a sequence generated by Algorithm INDO. Then 
 \begin{eqnarray}  \label{76}
 &\, & \|\uu^{k+1}-\uu^*\|^2_{\cal G} - \|\uu^{k} - \uu^*\|^2_{\cal G} \\
 \nonumber
& \leq &- \frac{2 \alpha}{m+M}\|\nabla f({\mx}^{k+1})-\nabla f({\mx}^*)\|^2 \\
&-& \|{\mx}^{k+1}-{\mx}^*\|^2_{(\frac{2 \alpha m M}{m+M} -\frac{\alpha}{\zeta})\II+\alpha^2(\II-\WW)} \nonumber \\
& -&  \alpha \varepsilon \|\md^k\|^2 +  \alpha \zeta \|\me^k\|^2. \nonumber
\end{eqnarray}
\end{lemma}

The main convergence result is stated in the following Theorem. 

\begin{te} \label{T} Assume that A\ref{A1}-A\ref{A3} hold and let $ \{{\mx}^k,\mq^k\} $ be a sequence generated by Algorithm INDO. Let $ \beta, \phi > 1 $ be arbitrary constants, $ \lambda_2 $ the smallest positive eigenvalue of $ (\II-\WW) $ and $ \zeta \in ((m+M)/(2mM),\varepsilon/(8M^2)).$ Then there exists $ \bar{\eta} > 0 $ such that for $ \eta_k \leq \bar{\eta}, $ the sequence of Lyapunov functions $ \|\uu^k-\uu^*\|_{\cal G} $ satisfies  
$$ \|\uu^{k+1}-\uu^*\|^2_{\cal G} \leq \frac{1+\tilde{\delta}}{1+\delta} \|\uu^k-\uu^*\|^2_{\cal G}, $$
where $\delta\leq \min\{\delta_{\mathrm{a}}, \delta_{\mathrm{b}}\}$ with  
\begin{eqnarray*} 
&\,& \delta_{\mathrm{a}}  =\frac{2mM}{(m+M)\varepsilon} - \frac{1}{\varepsilon \zeta}\\
&\,&
\delta_{\mathrm{b}} = \frac{(\alpha \varepsilon-8M^2 \alpha \zeta)(\phi-1)(\beta-1)\lambda_2}{\beta \varepsilon^2(\phi-1)+8M^2\beta(\beta-1)\phi},\frac{2 \alpha \lambda_2}{(m+M) \phi \beta},
\end{eqnarray*}
and $ \tilde{\delta} < \delta. $ 
\end{te}

{\em Proof. } To prove the statement we have to find  $ \delta $ and $ \tilde{\delta} $ such that 
\be \label{77} 
\delta \|\uu^{k+1} - \uu^*\|^2_{\cal G} - \tilde{\delta} \|\uu^k - \uu^*\|^2_{\cal G} \leq \|\uu^k - \uu^*\|^2_{\cal G} - \|\uu^{k+1} - \uu^*\|^2_{\cal G}. 
\ee
Using the estimate for $ \|{\mv}^{k+1}-{\mv}^*\| $ given in Lemma \ref{L4} and (\ref{76}) in Lemma \ref{L5}, the inequality (\ref{77}) holds if 
\begin{eqnarray}  \label{78} 
&\, & \delta \alpha \varepsilon \|{\mx}^{k+1}-{\mx}^*\|^2 + \frac{\delta \beta \varepsilon^2}{(\beta-1) \lambda_2} \|\md^k\|^2 \\
&+&
\nonumber
\frac{\delta \phi \beta}{\lambda_2} \|\nabla f({\mx}^{k+1}) - \nabla f({\mx}^*)\|^2   \\
& + & 
\nonumber \frac{\beta \phi \delta}{(\phi-1)\lambda_2}\|\me^k\|^2 - \tilde{\delta} \alpha \varepsilon \|{\mx}^k - {\mx}^*\|^2 - \tilde{\delta}\|{\mv}^k - {\mv}^*\|^2 
\\
& \leq &  \|{\mx}^{k+1}-{\mx}^*\|^2_{(\frac{2 \alpha m M}{m+M} -\frac{\alpha}{\zeta})\II+\alpha^2(\II-\WW)} + \alpha \varepsilon \|\md^k\|^2  \nonumber \\ & + & \frac{2 \alpha}{m+M} \|\nabla f({\mx}^{k+1})-\nabla f({\mx}^*)\|^2 - \alpha \zeta \|\me^k\|^2.  \nonumber
\end{eqnarray}

By Lemma \ref{L3} we have
\be \nonumber
\|\me^k\|^2  \leq  8 M^2\|\md^k\|^2 + 4 \eta_k^2 (M+ 2 \alpha)^2\|{\mx}^k-{\mx}^*\|^2 + 8 \eta_k^2 \|{\mv}^k - {\mv}^*\|^2. 
\ee
Substituting this bound for $ \|\me^k\|^2 $ and the corresponding lower bound for  $ -\|\me^k\|^2 $ at both sides of (\ref{78}) we get the inequality 
{\allowdisplaybreaks{
\begin{eqnarray} 
\label{79} 
0 & \leq & \|{\mx}^{k+1}-{\mx}^*\|^2_{(\frac{2 \alpha m M}{m+M} -\frac{\alpha}{\zeta})\II+\alpha^2(\II-\WW)} \\
&+&\nonumber  \|\nabla f({\mx}^{k+1})-\nabla f({\mx}^*)\|^2\left(\frac{2 \alpha}{m+M} - \frac{\delta \phi \beta}{\lambda_2}\right) \nonumber \\
& + & \|\md^k\|^2 \left(\alpha \varepsilon - \frac{\delta \beta \varepsilon^2}{(\beta-1) \lambda_2} - 8M^2 (\frac{\beta \phi \delta}{(\phi-1)\lambda_2} + \alpha \zeta)\right) \nonumber \\
& + & \|{\mx}^k - {\mx}^*\|^2\,(\,-4\eta^2_k(M+2\alpha)^2\left( \frac{\delta \beta \phi }{(\phi-1) \lambda_2} + \alpha \zeta\right) 
\nonumber \\
&+& \nonumber   \alpha \varepsilon \tilde{\delta} \,\,) \nonumber \\
& + & \|{\mv}^k - {\mv}^*\|^2\left(\tilde{\delta} - 8 \eta_k^2\left(\alpha \zeta + \frac{\delta \beta \phi }{(\phi-1) \lambda_2}\right)\right). \nonumber
\end{eqnarray}}}

The last inequality holds for $ \delta $ specified in the statement  and $  \tilde{\delta} < \delta $ if 
 \begin{eqnarray*}
 &\,& \eta^2_k \leq \bar{\eta}^2 \\
 &\,& \leq \min\{\frac{\tilde{\delta} \alpha \varepsilon}{4(M+2\alpha)^2(\alpha \zeta + \frac{\delta \beta \phi}{(\phi-1) \lambda_2})}, \frac{\tilde{\delta}}{8(\alpha \zeta + \frac{\delta \beta \phi}{(\phi-1) \lambda_2})}\} 
 \end{eqnarray*}
 $ \Box$

\section{Distributed computation  of Inexact Newton step } 

In this section we focus on the Step S1 in \eqref{58} of the algorithm INDO, for a fixed iteration $ k. $ We drop the iteration counter in the Hessian matrix $ \HH^k $ in (\ref{55}) and remaining relevant quantities to simplify notation.  That is, further on denote $ \HH^k = \HH, $ i.e., 
\be \label{81}
\HH = \nabla f({\mx}^k) + \alpha (\II - \WW) + \varepsilon \II, 
\ee
where $ \varepsilon >0 $ is the proximal parameter and $\mg=\mg^k=\nabla f(\mx^k) + \mq^k +  \alpha (\II - \WW) \mx^{k}$. By assumption A\ref{A2} the Hessian of the objective function $ f $ is positive definite with $ m\II \preceq \nabla^2 f({\mx}^k), $ the matrix   $  \alpha (\II - \WW) $ is positive semidefinite   and  since  $ \varepsilon>0 $ we conclude that there holds 
\be \nonumber (m+\varepsilon) \II \preceq \HH.
\ee  
Notice that  we can even make $\HH$  strictly diagonally dominant by taking $\varepsilon$  large enough.  Furthermore, $ \HH $ has the same sparsity structure of the network and hence one can easily apply an iterative solver of the fixed point type in S1 of Algorithm INDO. 

Let us look more closely at the step calculation in  S1. For simplicity of exposition here we concentrate on  the Jacobi Overrelaxation (JOR) method although other options for linear solver are possible. For each $ k $ we need to solve  the system 
$ \HH \md = -\mg $ approximately i.e. we are looking for $ \md^k $ such that (\ref{58}) holds.  
 Thus we will consider the linear system 
\be \HH \md = - \mg, \label{82} \ee
with $ \HH $ satisfying \eqref{81}. With the splitting 
$\HH = \DD - \GG, $ where $ \DD $ is the diagonal part of $ \HH, $ the Jacobi Overrelaxation matrix $ {\TT}_\gamma $ is defined as 
\be \label{83}
\TT_\gamma = \gamma \DD^{-1} \GG + (1-\gamma)\II 
\ee
and the JOR iterative method is defined as 
\be \label{84} \md^{\ell+1} = \TT_\gamma \md^\ell - \gamma \DD^{-1} \mg, 
\ee
with $ \gamma  $ being the relaxation parameter, for arbitrary $ \md^0. $ For $ \gamma = 1 $ we get the Jacobi iterative method with the iterative matrix  $ \TT_1 = \DD^{-1} \GG. $ 

The matrix $ \TT_{\gamma}$ has the same sparsity structure as $ \WW $ and all other matrices we considered so far, including the global Hessian. More precisely we can specify block-rows of $ \TT_\gamma $ that are used in the JOR method by each node as follows. Notice that, for matrix $DD$, we have its diagonal blocks $ D_{ii} \in \mathbb{R}^{n \times n}, i=1,\ldots,N, $ given by
\be \label{diagD} 
[D_{ii}]_j = \varepsilon + \alpha(1-w_{ii}) + [\nabla^2 f_i({\mx}_i^k)]_{jj}, j=1,\ldots,n, 
\ee
where  $[\nabla^2 f_i(x_i^k)]_{jj}$ denotes the $j$-th diagonal element of the local Hessian $ \nabla^2 f_i(x_i^k). $ Similarly, with matrix $\\G$, consider its blocks $ G_{ij} \in \mathbb{R}^{n \times n}, \; i,j=1,\ldots, N $. Then for $ i \neq j $ we have that $ G_{ij} $ is diagonal with elements $ \alpha w_{ij} \,I$ on the diagonal, and for $ i = j $ we have 
\be \label{diagG}
 G_{ii} = diag(\nabla^2 f_i(x_i^k))-\nabla^2 f_i(x_i^k)  , 
 \ee
  with $ diag(\nabla^2 f_i(x_i^k)) \in \mathbb{R}^{ n \times n}$ being the diagonal of local Hessian $ \nabla^2 f_i(x_i^k)), $ for all $ i=1,\ldots, N. $
Thus the matrix $ \TT_\gamma $ has block elements  $ [\TT_\gamma]_{ij} = D_{ii}^{-1} G_{ij}, \; i,j=1,\ldots,N. $  
Let $ [\TT_\gamma]_i \in \mathbb{R}^{n \times nN} $ be the block row of $ \TT_\gamma, \; i=1,\ldots,N$
\be \label{Mgamma}
[\TT_\gamma]_i = [D_{ii}^{-1} G_{i1},\ldots,D_{ii}^{-1} G_{ii},\ldots,D_{ii}^{-1} G_{iN} ]. 
\ee
Similarly, we can define the partition of the gradient vector $ \mg $ with the $ i$th component being 
\be \label{g} 
 g_i^k = \nabla f_i(x_i^k) + q_i^k + \alpha[(1-w_{ii})x_i^{k} - \sum_{j \in O_i} w_{ij} x_j^k], \;  i= 1,\ldots, N. 
 \ee
 Notice that each node $ i $ can compute $ [\TT_\gamma]_i$ and $ g_i^k.  $  
 
 Let us now state the $ k+1$-th iteration of INDO algorithm in the node-wise manner. 
 
 \noindent{\bf Algorithm INDO-nodewise}
\begin{itemize}
\item[Given]: $ {\mx}^k,\mq^k,  \alpha > 0, \; \varepsilon > 0, \ell_k \geq1, \ell \in \mathbb{N}.  $
\item[S1] Computing the direction $ \md^k $ 
\begin{itemize}
\item[S1.1] Each node computes 
$ [\TT_\gamma]_i $ and $ g_i^k $ by (\ref{diagD}) - (\ref{g}) and chooses $ d_i^0 \in \mathbb{R}^{n}. $ 
\item[S1.2] Each node computes $ d_i^k $  as follows. \\
For $\ell=0,\ldots, \ell_k-1$ \\
\hspace*{5pt} Each node sends $ d_i^{\ell} $ to all its neighbors and receives $ d_j^{\ell}, j\in O_i. $ \\
\hspace*{5pt} Each node computes  
\be \label{JORit}
d_i^{\ell+1} = [\TT_\gamma]_i \md^{\ell} - \gamma D_{ii}^{-1} g_i^k 
\ee
\hspace*{5pt} with $ \md^{\ell} = (d_1^{\ell},\ldots,d_N^{\ell}). $\\
Endfor.\\
Set $ d_i^k = d_i^{\ell_k}. $
\end{itemize}
\item[S2]  Primal update
\begin{itemize}
\item[S2.1]
Each node updates the primal variable
$$ x_i^{k+1} = x_i^k + d_i^k. $$
\item[S2.2]
Each node sends $ x_{i}^{k+1} $ to all its neighbours and receives $ x_j^{k+1}, j \in O_i. $ 
\end{itemize}
\item[S3]  
Each node computes  the dual update
$$ q_i^{k+1} = q_i^k + \alpha[(1-w_{ii}) x_i^{k+1} - \sum_{j \in O_i} w_{ij}x_j^{k+1}],$$ 
 and sets  $ k = k+1.   $ 
\end{itemize}

Notice that Step 1.2 requires only the neighboring elements of $ \md^{\ell} $ as $ G_{ij} = 0 $ if $ w_{ij} = 0, i\neq j,$ and thus the step is well defined and (\ref{JORit}) is equivalent to (\ref{84}). 
The algorithm above specifies a choice for Step S1 in the general INDO method in \eqref{58}, but it is not straightforward to understand the connection.  Namely, in Step 1.2 of the node-wise algorithm we state that each node should perform $ \ell_k $ JOR iterations, while the main algorithm (see \eqref{58}) requires the step $ \md^k $ such that the inexact forcing condition (\ref{58}) with some $ \eta_k > 0 $ holds. In the sequel we first analyse the convergence conditions of the JOR method for solving (\ref{82}) and then we show that one can in fact determine $ \ell_k $ such that (\ref{58}) holds after $ \ell_k $ iterations in Step 1.2 of the above algorithm. 

The JOR method is convergent for 
\be  \gamma \in (0, 2/\sigma(\DD^{-1}\HH)) \label{gammacond} \ee 
for symmetric positive definite matrices, with $ \sigma(\DD^{-1}\HH) $ being the spectral radius of  $ \DD^{-1}\HH, $ \cite{AG}. In the statement below we  estimate the upper bound for the spectral radius of matrix $\DD^{-1}\HH$, using the block-wise Euclidean norm as in \cite{DFIX, DQN}. For $ \TT \in \mathbb{R}^{nN\times nN}, \TT=[T_{ij}], T_{ij} \in \mathbb{R}^{n\times n} $ we define
\be \label{blocknorm}
\|\TT\|_{\mathrm{b}} =\max_{1\leq i\leq N}  \sum_{j=1}^N \|T_{ij}\|_2. 
\ee
Clearly, $ \|\cdot\|_{\mathrm{b}} $ is a norm.

\begin{prop} \label{Pgamma}
Let Assumptions \ref{A1}-\ref{A3} hold and $ \alpha > 0. $ Then the Jacobi Overrrelaxation method (\ref{82}) is convergent for 
$$ \gamma \in \left(0, 2 \frac{m +\alpha(1-w_d) + \varepsilon }{M+ \varepsilon + \alpha(1-w_m) + \alpha(1-w_d)}\right), $$
where $ w_d = \max_{1\leq i \leq N} w_{ii}, \; w_m = \min_{1\leq i \leq N} w_{ii}. $
\end{prop}

{\em Proof.}  Given the fact that for positive $ \varepsilon $ and $ \alpha$ we have that $ \HH $ is positive definite, the convergence interval is determined by (\ref{gammacond}). On the other hand we have $ \sigma(\DD^{-1} \HH) \leq \|\DD^{-1} \HH\|_{\mathrm{b}}, $ so we have to estimate $ \|\DD^{-1} \HH\|_{\mathrm{b}}. $   By definition of $ \DD^{-1} \HH $ and the norm we have
\begin{eqnarray*}  &\,&  \|\DD^{-1} \HH\|_{\mathrm{b}}  =  \max_{1\leq i\leq N} \sum_{j=1}^N \|[\DD^{-1} \HH]_{ij}\|_2 
\\ & = & \max_{1\leq i\leq N} \left(\|I - D_{ii}^{-1} G_{ii}\|_2 + \sum_{j\neq i} \|-D_{ii}^{-1} G_{ij}\|_2 \right). 
\end{eqnarray*}
For each $ i=1,\ldots, N $ we have that $ D_{ii} $ is given by (\ref{diagD}). Furthermore, the diagonal elements of local Hessian $ \nabla^2 f_i(x_i^k) $ are in the interval $ [m,M] $ by A\ref{A2}. Therefore 
$$ \|D_{ii}\|_2 \geq m + \alpha (1-w_{ii}) + \varepsilon \geq m + \alpha (1-w_d) + \varepsilon. $$ 
Next, for $i \neq j$, 
   we have $\|G_{ij}\|_2 = \alpha\,w_{ij}$ and $ \sum_{j \in O_i} w_{ij} = 1- w_{ii}  \leq 1- w_m $ by the properties of $ W $ stated in A\ref{A1}. 
   On the other hand, for 
   $i=j$, we have by (\ref{diagG})  
   $$\|G_{ii}\|_2
   \leq 
 \|\nabla^2 f_i - diag(\nabla^2 f_i)\|_2.$$
   It can be shown that 
   $\|\nabla^2 f_i(x_i^k) - diag(\nabla^2 f_i(x_i^k))\|_2
    \leq M-m$ and thus 
$$\|G_{ii}\|_2 \leq M - m. $$
Combining the bounds for $ D_{ii} $ and $ G_{ii} $    we get 
\begin{eqnarray*} 
\|I - D_{ii}^{-1} G_{ii}\|_2 & \leq & 1 + \|D_{ii}^{-1} G_{ii}\|_2  \\
&\leq& 1 + \frac{M-m}{m + \alpha(1-w_d) + \varepsilon} \\
& = & \frac{M + \alpha(1-w_d) + \varepsilon}{m + \alpha(1-w_d) + \varepsilon}.
\end{eqnarray*} 
Finally, using the bound above and the bounds for $ \|D_{ii}^{-1}\|_2 $ and $ \|G_{ij}\|_2, $ we get
\begin{eqnarray*} 
 \|\DD^{-1} \HH\|_{\mathrm{b}} & \leq & \frac{\alpha(1-w_m)}{m + \alpha(1-w_d)+\varepsilon} +\frac{M + \alpha(1-w_d) + \varepsilon}{m + \alpha(1-w_d) + \varepsilon} \\
 & = & \frac{M + \varepsilon + \alpha(1-w_m) + \alpha(1-w_d)}{m + \alpha(1-w_d)+\varepsilon},
 \end{eqnarray*}
 and the statement follows. $ \Box $
  
  Notice that the bound obtained above does not depend on $ k, $ and hence we can claim that the method in Step 1.2 of the node-wise algorithm will converge to the true solution $ \md^* $ for an arbitrary $ \md^0 \in \mathbb{R}^{nN} $ if we define the JOR matrix with the relaxation parameter $ \gamma $  as stated in Proposition \ref{Pgamma}. 
  
  The remaining question is the number of inner iterations in Step 1.2 $ \ell_k $ that ensures that the forcing condition (\ref{58}) is satisfied. In the statement bellow we will prove that such $ \ell_k $ exists considering a special choice of $ d_i^0 $ for simplicity of exposition. Similar result can be proved for other choices of $ d_i^0. $ 
  
  \begin{prop} \label{Pelk}
  Let  A\ref{A1}-A\ref{A3} hold and assume that $ \gamma $ is chosen from the interval specified in Proposition \ref{Pgamma}. For given $ \eta_k > 0 $ denote by $ \md^k $ a step  obtained through Step 1.2 of the node-wise algorithm INDO with $ d_i^0 = 0, i=1,\ldots, N  $ and $ \ell_k $ inner iterations (\ref{JORit}).  Then $ \|\HH \md^k + \mg\| \leq \eta_k \|\mg\| $ after at most 
  $$ \ell_k \geq  \left \lceil \frac{\ln (\eta_k /c)}{|\ln
(\sigma(\TT_\gamma))|}\right \rceil,$$ where 
  iterations, where $\lceil a \rceil$ 
  denotes the smallest integer greater than or equal to $a$ and $$ c:=\frac{M+2+\varepsilon}{m+\varepsilon} \sqrt{\frac{\varepsilon+M+\alpha (1-w_m)}{\varepsilon+m+\alpha
(1-w_d)}}.$$   
  \end{prop} 
  
  {\em Proof. }
  Denote by $ \md^* $ the exact solution of $ \HH \md = - \mg $ and consider the sequence $ \{d^\ell\} $ generated at Step 1.2 with $ \md^0 = 0. $ As $ \md^{\ell+1} = \TT_{\gamma} \md^{\ell}-\gamma \DD^{-1} \mg, \; \ell=0,1,\ldots $  and $\md^*=\TT_{\gamma} \md^*-\gamma \DD^{-1} \mg$ we have 
  $$\md^{\ell} - \md^*=\TT_{\gamma}^{\ell}(\md^0 - \md^*).$$
  Notice that $\TT_{\gamma}$ is not symmetric, but it has the same set of eigenvalues as the following symmetric matrix  $$\TT'_{\gamma}=\DD^{-1/2} (\gamma \GG+(1-\gamma)\DD) \DD^{-1/2}:=\DD^{-1/2} \CC_{\gamma} \DD^{-1/2},$$
  since $\TT_{\gamma}=\DD^{-1} \CC_{\gamma}$ and both $\DD^{-1}$ and $\CC_{\gamma}$ are symmetric (see Remark 4.2 of \cite{koreni} for instance). Therefore, we conclude 
  $$\sigma(\TT_{\gamma})=\sigma(\TT'_{\gamma})=\|\TT'_{\gamma}\|$$
  and according to the choice of $\gamma$ and Proposition \ref{Pgamma} we have $\sigma(\TT_{\gamma})<1$. Moreover, it ca be shown that $$\TT_{\gamma}^{\ell}=\DD^{-1/2} (\TT'_{\gamma})^{\ell} \DD^{1/2}$$
  and thus we obtain 
  \begin{eqnarray} \label{dl} 
   &\,&\|\md^{\ell} - \md^*\| \leq \|\DD^{-1/2}\| \|\TT'_\gamma\|^{\ell} \|\DD^{1/2}\| \|\md^0 - \md^*\| \\
   &= &  \|\DD^{-1/2}\| (\sigma(\TT_{\gamma}))^{\ell} \|\DD^{1/2}\| \|\md^0 - \md^*\|.  \nonumber 
   \end{eqnarray}
   Considering that, by the structure of $\DD$, we have 
   $$(\varepsilon+m+\alpha (1-w_d))\II \preceq \DD \preceq (\varepsilon+M+\alpha (1-w_m))\II,$$
   we conclude that 
   $$\|\DD^{-1/2}\| \|\DD^{1/2}\| \leq \sqrt{\frac{\varepsilon+M+\alpha (1-w_m)}{\varepsilon+m+\alpha
(1-w_d)}}:=c_D.$$ 
Therefore, we obtain 
 \begin{eqnarray} \label{dl2} 
   &\,&\|\md^{\ell} - \md^*\| \leq  (\sigma(\TT_{\gamma}))^{\ell}  c_D  \|\md^0 - \md^*\| \\
   &= & (\sigma(\TT_{\gamma}))^{\ell}  c_D  \|\HH^{-1} \mg \| \\
   &\leq  & (\sigma(\TT_{\gamma}))^{\ell}  c_D  \|\HH^{-1}\| \| \mg \| \\
   &\leq  & (\sigma(\TT_{\gamma}))^{\ell}   \frac{c_D }{m + \varepsilon } \| \mg \| , 
   \end{eqnarray}
   since by A\ref{A1} and A\ref{A2} we can upper bound the Hessian and its inverse as 
   $$ \|\HH\| \leq M + 2 +\varepsilon, \; \|\HH^{-1}\| < (m+\varepsilon)^{-1}, $$ as $ \|\II - \WW\| \leq 2.$ Finally,  we have
   \begin{eqnarray*}
    &\,&\|\HH \md^{\ell} + \mg\| = \|\HH \md^{\ell} - \HH \md^*\| \leq \|\HH\|\|\md^{\ell} \\
    &-& \md^*\| \leq c_D \frac{M+2 + \varepsilon}{m+\varepsilon} (\sigma(\TT_\gamma))^{\ell}\|\mg\| 
    \end{eqnarray*}
    and we get 
   $$ \|\HH \md^{\ell_k} + \mg\| \leq \eta_k \|\mg\| $$ for $ \ell_k $ given in the statement. $ \Box $
   
   The above statement implies that we can run INDO with the fixed number of inner iterations for all $k, $ by taking $ \ell_k = \ell(\bar{\eta})  $ (see Theorem~\ref{T}),   to avoid the overhead needed to define suitable 
   iteration-varying quantities $ \eta_k $ and $ \ell_k $.  In fact, the numerical results we present in the next Section show that INDO behaves remarkably well if we take a very modest number of inner iterations at each outer iteration; namely, even $ \ell_k =1 $ with the so-called warm start gives satisfactory results. 
   In practice, this also alleviates the need to calculate the complicated expression for $\ell(\bar{\eta})$ 
   as per Proposition~{4.2} that depends on certain global quantities such as $\sigma(\TT_\gamma)$.

   Another possibility  would be to define the forcing condition in the block infinity norm, i.e.,  
   $$ \|\HH \md^k + \mg\|_{\infty} \leq \eta_k \|\mg\|_{\infty}. $$ Then each node can compute the value of local residual $ \|r_i^{\ell}\|_{\infty} $ at each $\ell$; subsequently, at each $\ell$, the nodes can run an iterative primitive for computing  maximum (e.g., \cite{Max}) of their local residuals until the overall forcing condition in the block infinity norm is satisfied.
   
   A couple of words on the forcing sequence $ \{\eta_k\} $ are due here. Given that we are in the framework of the proximal method of multipliers, the most we can achieve is the overall linear convergence as the linear convergence is achieved even if the Newtonian linear systems are solved exactly i.e. for $ \eta_k = 0, $ see \cite{ESOM}.  On the other hand we know from the classical optimization theory that $ \eta_k $ greatly influences the local convergence rate of Inexact Newton methods. In fact for $ \eta_k \to 0 $ we have superlinear local convergence while $ \eta_k = {\cal O}(\|\mg^k\|) $ recovers the quadratic convergence of the Newton method. In fact, one can show that the bound for $ \bar{\eta} $ can be made arbitrarily close to 1 and still have a convergent method if a weighted norm is used, see \cite{ypma}. But given the overall linear convergence at most in the assumed framework, it seems most efficient to work with a constant $ \eta_k $ during the whole process as such procedure maintains the linear rate of convergence while minimizing the effort needed for tuning $ \eta_k $ and $ \ell_k. $ The key property -- avoiding oversolving in the approximate Newton step calculation -- is achieved with a constant $ \eta_k $ (and hence constant $ \ell_k).$ The comparison with ESOM, both in terms of (inner iterations) convergence factor and numerically, presented in the next section strongly supports the approach we advocate here -- approximate second order direction with the error proportional to the gradient value while incurring a small computational cost.

\section{Performance analysis: Numerical results and further analytical insights}

In this Section we analyze computational and communication costs of INDO and test its performance on a set of standard test examples. INDO method is a second-order method within Proximal Method of Multipliers that differs from the ESOM method, \cite{ESOM} in computing the primal step. So, we will first compare the theoretical convergence rates of these two methods.  These findings are presented in Subsection 5.1, while the numerical examples and cost comparison are presented in Subsection 5.2. 

\subsection{Inner iterations' convergence rate of INDO and ESOM}

We provide here a more detailed comparison of INDO with the ESOM method proposed in \cite{ESOM} with respect to inner iterations' convergence. The two algorithms have identical dual variable updates, given by \eqref{dualq}. Hence we focus on the primal variable updates.
 At each outer iteration~$k$, both INDO and ESOM approximately solve the system of linear equations~\eqref{82}. 
 While INDO solves~\eqref{82} via 
  \eqref{83}--\eqref{84}, 
  \cite{ESOM} adopts a different approach through a Taylor approximation; see equation~(13) in \cite{ESOM}. However, the Taylor approximation approach admits a representation in the spirit of \eqref{83}--\eqref{84}. Namely, it can be shown that 
  the ESOM solver of~\eqref{82} can be expressed as follows: %
\be \label{ESOM-inner-1} \md_{E}^{\ell+1} = \left( \DD_E^{-1} \BB\right) \md_E^\ell - \DD_E^{-1} \mg, 
\ee
for $\ell=-1,0,1,...,$ with 
$\md_E^{-1}=0$. 
Here, ESOM utilizes the splitting 
$\HH=\DD_E - \BB$ of the Hessian $\HH$ in \eqref{81}, 
with
\begin{eqnarray}
\DD_E & = & \nabla^2 f({\mx}^k) + 
2\alpha \,(I -\,\mathrm{diag}(\WW)) 
+ \varepsilon\,I\\
\BB &=& \alpha \,(I-2\mathrm{diag}(\WW) +\WW).
\end{eqnarray}
That is, an ESOM inner iteration $\ell$ in \eqref{ESOM-inner-1}  is of a form similar to INDO; however, a major difference is that the utilized splitting involves a block-diagonal matrix $\DD_E$ with non-sparse $n \times n$ diagonal blocks. As a result, \eqref{ESOM-inner-1} requires inverting a dense $n \times n$ positive definite matrix per node. In contrast, with the INDO approach in \eqref{83}--\eqref{84}, the splitting matrix $\DD$ is diagonal and hence no matrix inversion is required. A more detailed computational cost analysis per inner iteration for quadratic and logistic regression problems is presented in subsection 5.2. 

We next compare convergence rates of inner iterations of ESOM and INDO. For both methods, it is easy to show that the error $\mc^{\ell} = \md^\ell - \md^*$ with respect to the solution $ \md^* $ of~\eqref{81} evolves as:
\begin{equation}
    \mc^{\ell+1} = \TT \mc^{\ell}.
\end{equation}
Here, for 
ESOM, we have that $\mc^{\ell} = \md_E^{\ell}-\HH^{-1}\mg$,
 and $\TT = \DD_E^{-1} \BB$; and for INDO, 
 we have 
 $\mc^{\ell} = \md^{\ell}-\HH^{-1}\mg$,
 and $\TT = \TT_{\gamma} = \gamma\,\DD^{-1}\GG+(1-\gamma)\II$.
 That is, provided that the spectral radius  
 $\sigma(\TT)$ is less than one, $\mc^{\ell}$ 
 converges to zero linearly with the convergence factor determined by~$\sigma(\TT)$. It is in general difficult to explicitly evaluate $\sigma(\TT)$ for INDO and ESOM. We hence compare the two methods in terms of the block-wise matrix norm upper bound 
 $\|\TT\|_{\mathrm{b}} $ on $\sigma(\TT)$ introduced in (\ref{blocknorm}), that in turn admits 
 intuitive upper bounds.  
 Assuming for simplicity 
 that all elements on the diagonal of $W$
  are mutually equal, $w_{ii} = w,\, i=1,\ldots,N$, and 
 that $\gamma=1$ with INDO, it can be shown,  similarly as in Proposition \ref{Pgamma},  
   that the following holds:
 \begin{eqnarray}
 \label{eqn-M-b-norms-ESOM}
 & &\mathrm{ESOM}:\,\,
\|\TT\|_{\mathrm{b}} \leq 
\frac{2 \alpha \,(1-w)}
{2 \alpha \,(1-w)  + \varepsilon + m}<1\\
\label{eqn-M-b-norms-INDO}
& &
\mathrm{INDO}:  \|\TT\|_{\mathrm{b}} \leq 
\frac{M-m +  \alpha \,(1-w)}
{\alpha \, (1-w)  + \varepsilon + m}<1,\\
&\,&
\mathrm{for} \,\varepsilon>\max\{M-2m,0\}.
\nonumber
 \end{eqnarray}

We now comment on the convergence factor upper bounds in \eqref{eqn-M-b-norms-ESOM} and \eqref{eqn-M-b-norms-INDO}. First, we can see that, with both ESOM and INDO, the convergence factor can be made arbitrarily good (close to zero) by taking a sufficiently large $\varepsilon$. 
However, a too large $\varepsilon$ comes at a price of slowing down the outer iterations, i.e., as too large $\varepsilon$ makes small differences between ${\mx}^{k+1}$ and ${\mx}^k$, see equations (\ref{dl}) and (\ref{primal}). 
Second, we can see 
that, while we can take arbitrary (in fact, arbitrarily small) positive $\varepsilon$ for ESOM, with INDO $\varepsilon$ needs to be sufficiently large, i.e., we need to have $ \varepsilon>\max\{M-2m,0\}$. That is, 
by avoiding $n \times n$ matrix inversion at each node with INDO, we pay a price in that $\varepsilon$ should be sufficiently large, i.e., of order~$M$. 
However, extensive simulations in the next subsection show  that this incurs no   
overall loss of INDO, i.e., INDO is comparable or faster iteration and communication-wise (and faster computational cost-wise) than ESOM with a best hand-tuned~$\varepsilon$. Third, interestingly, when 
$\alpha$ and $\varepsilon$ are large compared with $M$, INDO's inner iteration convergence factor upper bound  
\eqref{eqn-M-b-norms-INDO} is comparable or even smaller than that of ESOM in \eqref{eqn-M-b-norms-INDO}. 
Extensive simulations show that 
taking $\alpha$, $\varepsilon$ to be of the same order and of the same order as $M$ (in fact, we can take $\alpha=M=\varepsilon$) works well with INDO. This choice is in a good agreement with the theoretical upper bound in~\eqref{eqn-M-b-norms-INDO}. 

The recommended tuning parameter choice $\alpha=M=\varepsilon$ still requires a beforehand global knowledge of system parameters, specifically the constant~$M$. Assuming that each node 
knows $M_i$--a Lipschitz constant of its own function $f_i$'s gradient, 
  $M=\max_{i=1.,,,,N}M_i$ can be obtained by running beforehand a distributed
 algorithm for computing maximum of scalar values held by the nodes, e.g.,  \cite{Max}.  
 This can be done with a low overhead, wherein the number of required  inter-neighbor (scalar-transmission) communication rounds is on the order of the network diameter. 
 INDO as implemented in Section~V also requires the beforehand knowledge of $m$ and $w_d$; these quantities can be computed beforehand analogously to~$M$.  

\subsection{Numerical results}

 We compare INDO with ESOM proposed in \cite{ESOM} on quadratic cost functions (simulated data) and on logistic regression problems (real data). 
The network we consider has $ N = 30 $ nodes and is formed as follows, \cite{Espectral}. The  points are sampled  randomly and uniformly from $[0,1] \times [0,1]$. The edges between points exists if their Euclidean distance is smaller than  $r=\sqrt{\log(N)/N}$. 
The resulting graph instance considered is connected. The weight coefficients in the communication matrix $ W $ are taken as  $w_{i,j}=1/(1+\max \{deg(i),deg(j)\})$, where $deg(i)$ stands for the degree of node $i,$ for directly connected nodes $ i $ and $ j, $ and the diagonal weights are $ w_{i,i}=1-\sum_{j\neq i}w_{i,j}$. The matrix $ W $ generated in this way satisfies A\ref{A1}. 

Let us now describe the test examples. 
Quadratic local cost functions are of the form 
$$f_i(y)=\frac{1}{2} (y-b_i)^T B_{ii} (y-b_i)$$
and the data is simulated as in   \cite{Espectral}, i.e.,  vectors $b_i$ are drawn from the Uniform distribution on $[1,31]$, independently from each other. Matrices $B_{ii}$ are of the form $B_{ii}=P_i S_i P_i$, where $S_i$ are diagonal matrices with Uniform distribution on $[1,101]$ and  $P_i$ are matrices of orthonormal eigenvectors of $\frac{1}{2} (C_i+C_i^T)$ where  $C_i$ have components drawn independently from the standard Normal distribution. Given that in this case we can compute the exact minimizer $ y^*, $ the error is measured as 
\be \label{errore} E({\mx}^k):=\frac{1}{N} \sum_{i=1}^{N}\frac{\|x_i^k-y^*\|}{\|y^*\|},\ee
Both Algorithms require the Hessian lower and upper bound $ M, m $ and we calculate them as  $M=\max_{i} M_i$, where $M_i$ is the largest eigenvalue of $B_{ii}$ and $m=\min_{i} m_i$, where $m_i>0$ is the smallest eigenvalue of $B_{ii}$. The dimension of the problem is set to  $n=100$.

The two methods are also compared on  binary classification problems and the following data sets: Mushrooms \cite{mush} ($n=112$, total sample size $T=8124$), LSVT Voice Rehabilitation \cite{VOICE} ($n=309$, total sample size $T=126$) and Parkinson's Disease Classification \cite{PARKINSON} ($n=754$, total sample size $T=756$). For each of the problems, the data is divided across 30 nodes of the graph described above.  The logistic regression with the quadratic regularization is used and thus the local objective functions are of the form 
$$f_i(y)=\frac{1}{|J_i|} \sum_{j \in J_i} \log(1+e^{-\zeta_j p_j^T y})+\frac{m}{2}\|y\|^2,$$
where $J_i$  collects the indices of the data points assigned to node  $i$, $p_j \in \RR^n $ is the corresponding vector of attributes and $\zeta_j \in \{-1,1\}$ represents the label.  
The data is scaled in a such way that $M=1+m$ with  $m=10^{-4}$. Since the solution is unknown in general, the error is measured as the average value of the original objective function  across the nodes' estimates  
\be \label{errorv} V(\mx^k)=\frac{1}{N} \sum_{i=1}^{N} \sum_{j=1}^{N} f_j(x_i^k).\ee 

We fix the free parameters of INDO method to $\alpha=\varepsilon=M$. The reasoning behind this choice is explained in Section 5.1. We use this choice in all the tested examples, although it may not be the optimal.  As discussed in the previous section, step S1 of INDO method can be implemented in different ways depending on the network characteristics, dimension of the problem, etc. We use a practical version of INDO in the tests by taking a fixed number $\ell$ of inner  iterations of JOR method (e.g., $\ell=1$)  with $\gamma=2(m+\varepsilon+\alpha(1-w_d))/(M+2\alpha+\varepsilon)$ denoted by  INDO-$\ell$ in the sequel.  Initial  $\md^{0}$ is obtained by solving $\|\HH^0 \md+\mg^0\|_{\infty} \leq \|\mg^0\|_{\infty}$. In all the subsequent iterations we use the so called warm start, i.e., we set $\md^{0}=\md^{k-1}$ at Step 1.2 in the node-wise INDO representation. Given that we implement the fixed number of inner iterations $ \ell$ the forcing term $ \eta_k $ is not explicitly imposed. 

Denote by ESOM-$\ell$-$\alpha$-$\varepsilon$ the ESOM method with $\ell$ inner iterations and the corresponding free parameters $\alpha$ and $\varepsilon$. Notice that ESOM-$\ell$ requires the same amount of communications as INDO-$\ell$. Therefore, we plot  iterations versus the corresponding error measures of both methods to compare the performance with respect to the communication costs. 

Regarding the computational cost, the main advantage of INDO method with respect to ESOM is the following: ESOM requires  inverting  full $ n \times n $ matrices at each node, while INDO only inverts the diagonal ones. The difference is more evident in non-quadratic case where the corresponding matrices need to be inverted in every outer iteration of ESOM method. This is even more significant in problems with relatively large $n.$  

Let us now estimate the computational costs of the two algorithms more precisely. Both ESOM and INDO are second order methods, so the cost of calculating the derivatives are the same. Moreover, the outer iterations updates are identical. The main difference lies in performing the inner iterations. In order to estimate the costs, let us observe the formula for INDO \eqref{JORit}, i.e., 
\be 
d_i^{\ell+1} = \gamma D^{-1}_{ii} \left( G_{ii} d_i^\ell  +\alpha \sum_{j \in O_i}   w_{ij} d_j^\ell - g_i^k\right)+(1-\gamma)d_i^\ell. 
\ee 
and ESOM
\begin{eqnarray} 
\nonumber 
&\,&d_i^{\ell+1} = [D_{E}^{-1}]_{ii} \left( \alpha (1-w_{ii}) d_i^\ell  +\alpha \sum_{j \in O_i}  w_{ij} d_j^\ell -  g_i^k\right),\\
&\,&\quad  d_i^{0} = -[D_{E}^{-1}]_{ii} g^k_i.
\label{ESOMd}
\end{eqnarray}
Both methods perform one consensus step ($\sum_{j \in O_i}   w_{ij} d_j^\ell$) per inner iteration per node. Both methods perform one  matrix-vector product per inner iteration per node as well: INDO calculates $G_{ii} d_i^\ell$ and ESOM calculates product of $[D_{E}^{-1}]_{ii}$ with the corresponding vector. The difference lies in the following. INDO inverts the diagonal matrix $D_{ii}$ and multiplies it with the corresponding vector in each inner iteration, i.e., at each inner iteration we have component-wise division of two vectors. Thus, the cost of this operations can be estimated to $n \ell$ scalar products in $\mathbb{R}^n$ (SPs) per outer iteration per node for INDO-$\ell$ algorithm. On the other hand, ESOM calculates the inverse of possibly dense symmetric positive definite matrix  $[D_{E}]_{ii}$ which can be estimated to $n^2 /6 $ SPs. For logistic regression problems, we estimate the common computational costs (see \cite{EFIX} for more details) of both algorithms to $|J_i|(2+n/2)+N+n \ell+N \ell/n$ SPs per node per (outer) iteration: $|J_i|(2+n/2)$ SPs is the cost of calculating the gradient and the Hessian of local cost function, $N=n N/n$ SPs comes from consensus step in calculating $\mg^k$, matrix-vector products take $n \ell$ SPs and the consensus with respect to $d_j$ vectors takes $\ell N/n$ SPs. Thus, the overall computational cost per node per iteration in logistic regression case is estimated to 
$|J_i|(2+n/2)+N+n\ell+N\ell/n+n\ell$
for INDO-$\ell$ and to 
$|J_i|(2+n/2)+N+n\ell+N\ell/n+n^2/6$
for ESOM-$\ell$. Clearly, these costs differ by the order $ {\cal O}(n). $ For the quadratic costs we do not have costs of calculating the local Hessian in every iteration while the cost of calculating the local gradient is $n$ SPs. Other common costs are the same as in logistic regression case. INDO still needs $n\ell$ SPs in every iteration, while ESOM only inverts the Hessian in the initial phase with the cost of $n^2/6$ SPs.   

\begin{figure}
  \includegraphics[width=5.5cm, height=9cm, angle=-90]{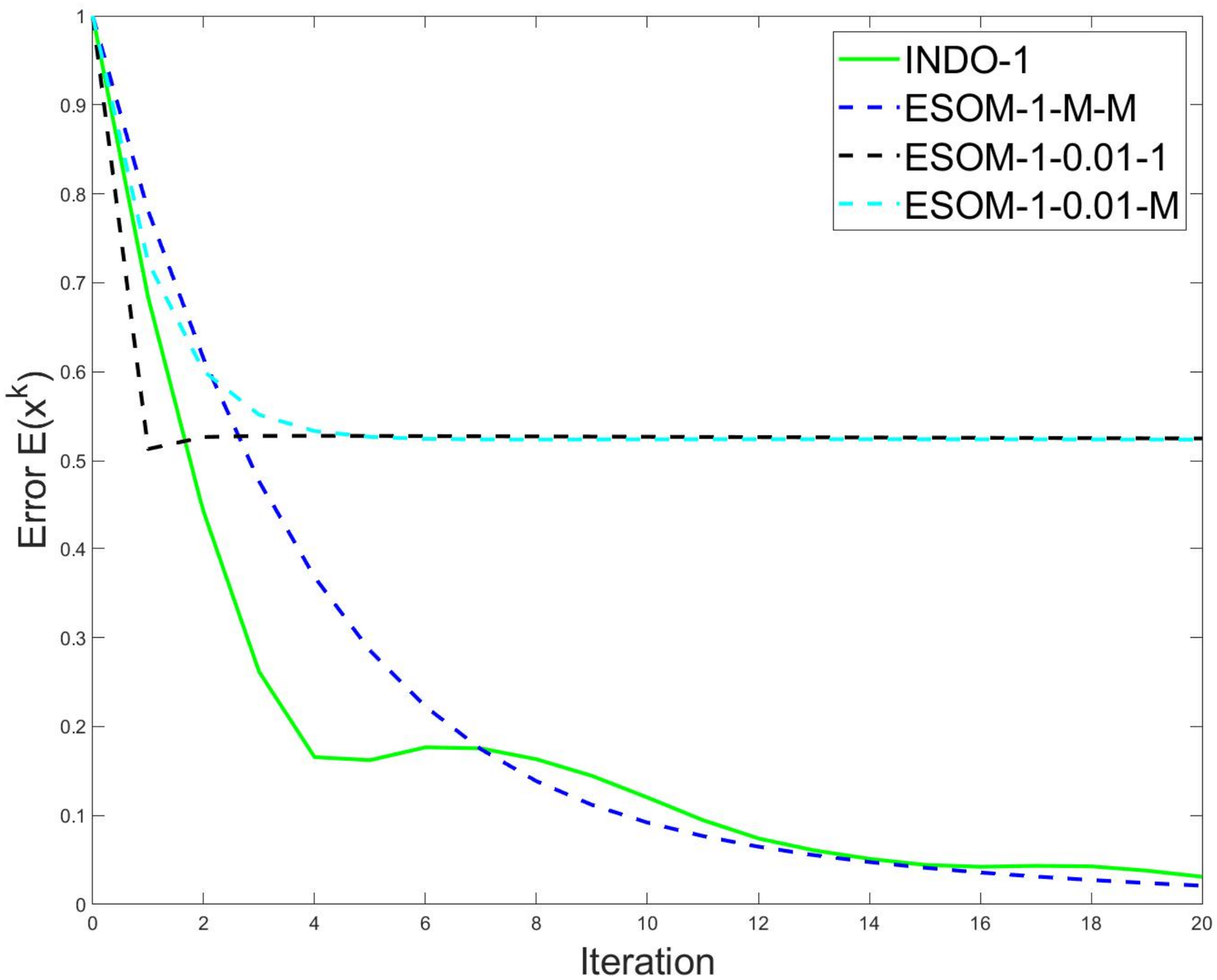}
  \includegraphics[width=5.5cm, height=9cm, angle=-90]{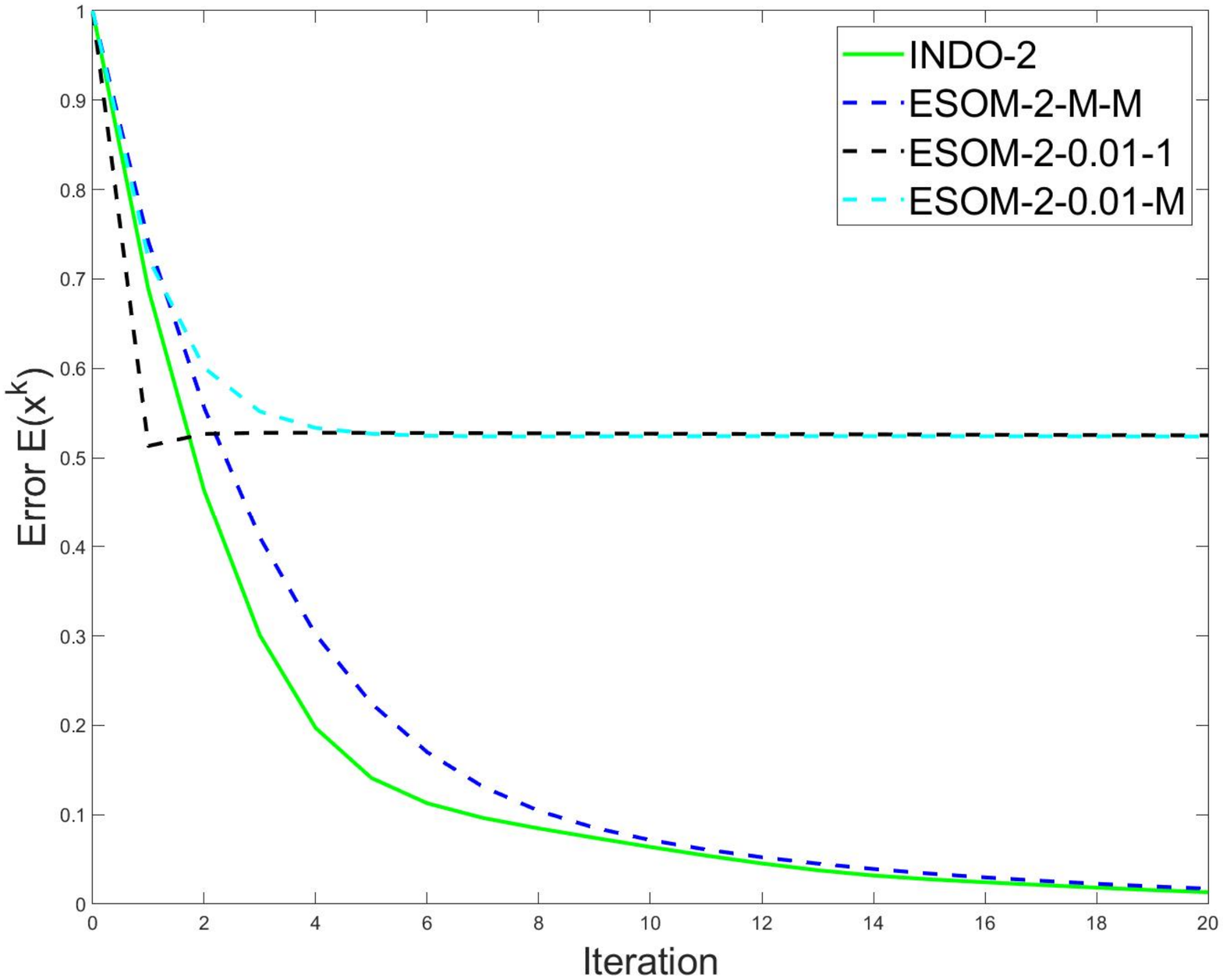}
\includegraphics[width=5.5cm, height=9cm, angle=-90]{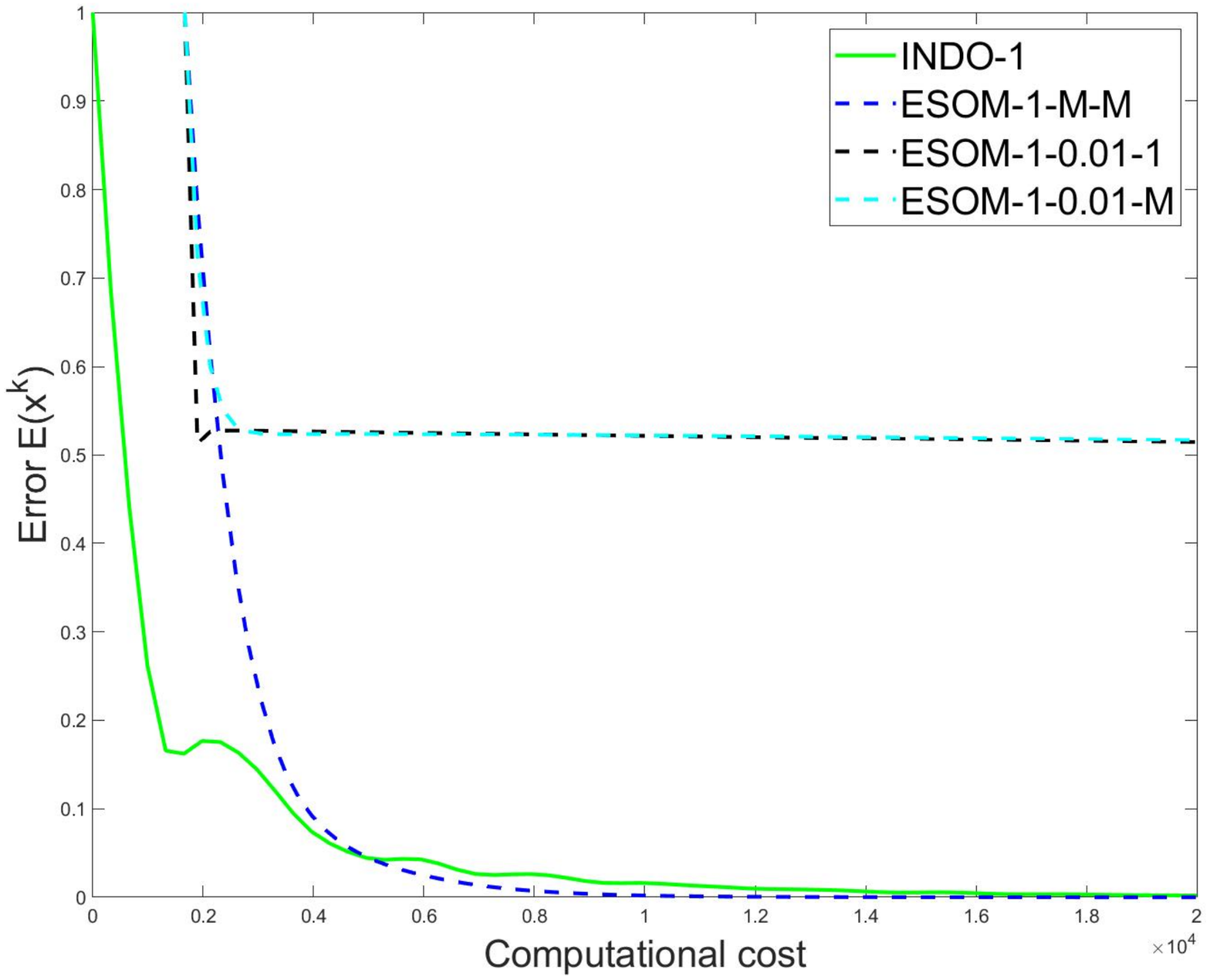}
 \includegraphics[width=5.5cm, height=9cm, angle=-90]{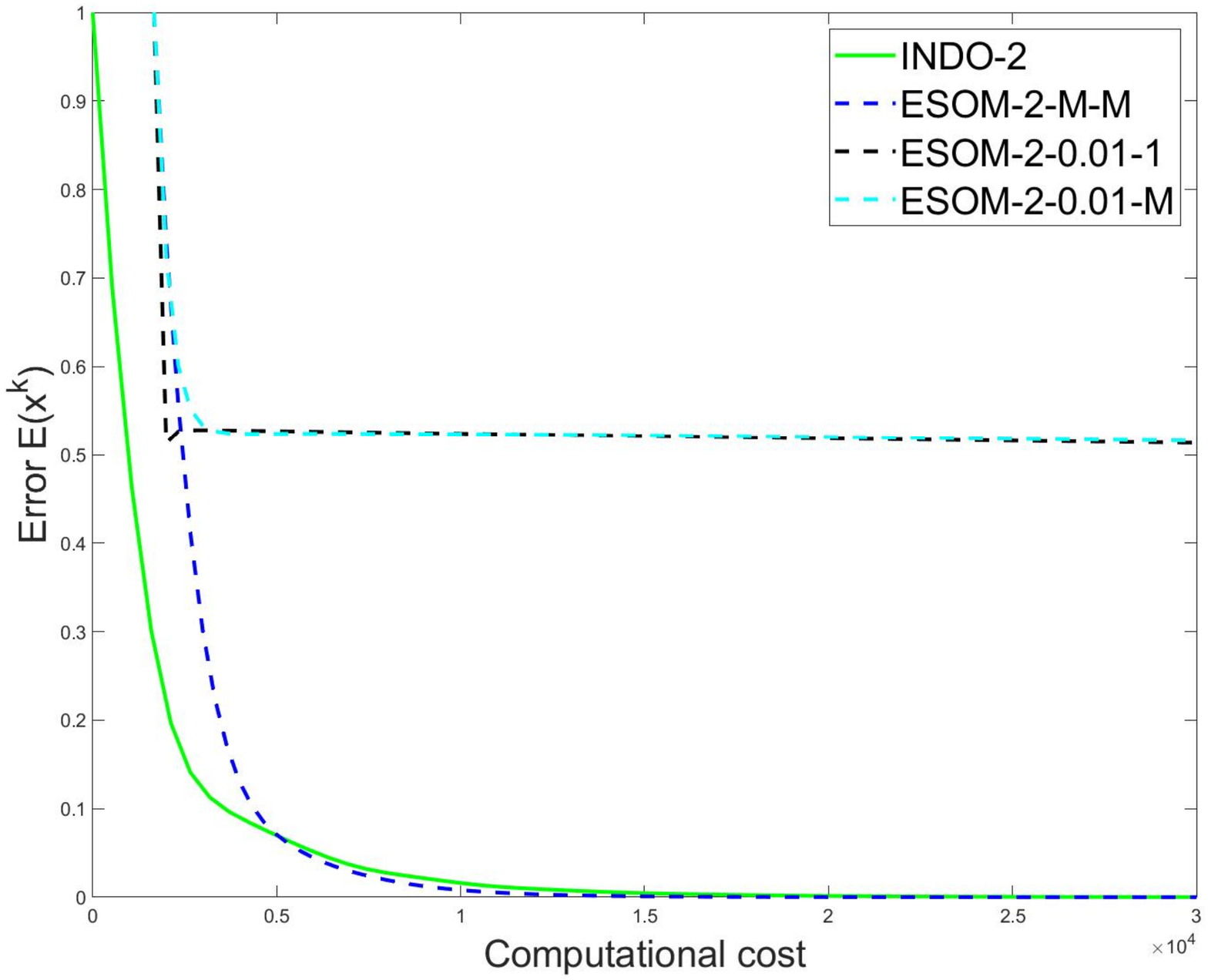}
 \caption{INDO (solid line) versus  ESOM methods (dotted lines): error \eqref{errore} with respect to iterations/communication cost (the two top figures) and computational cost (the two bottom figures) for the number of inner iterations $\ell=1$ (first and third figure from top) and  $\ell=2$ (second and fourth figure from top);  Simulated quadratic costs with  $n=100$ and  $N=30$.}
\end{figure}

Figure 1 presents results on simulated quadratic costs for different number of inner iterations, i.e. with respect to the communication costs. We test different combinations of $\alpha$ and $\varepsilon$ for ESOM method, while INDO is tested with fixed parameters as explained above.  
The best tested ESOM algorithm is  ESOM-$\ell$-M-M. Other variants of ESOM start better, but they  fail to converge to the same solution vicinity of the exact solution as the INDO method.  

\begin{figure}
\includegraphics[width=5.5cm, height=9cm, angle=-90]{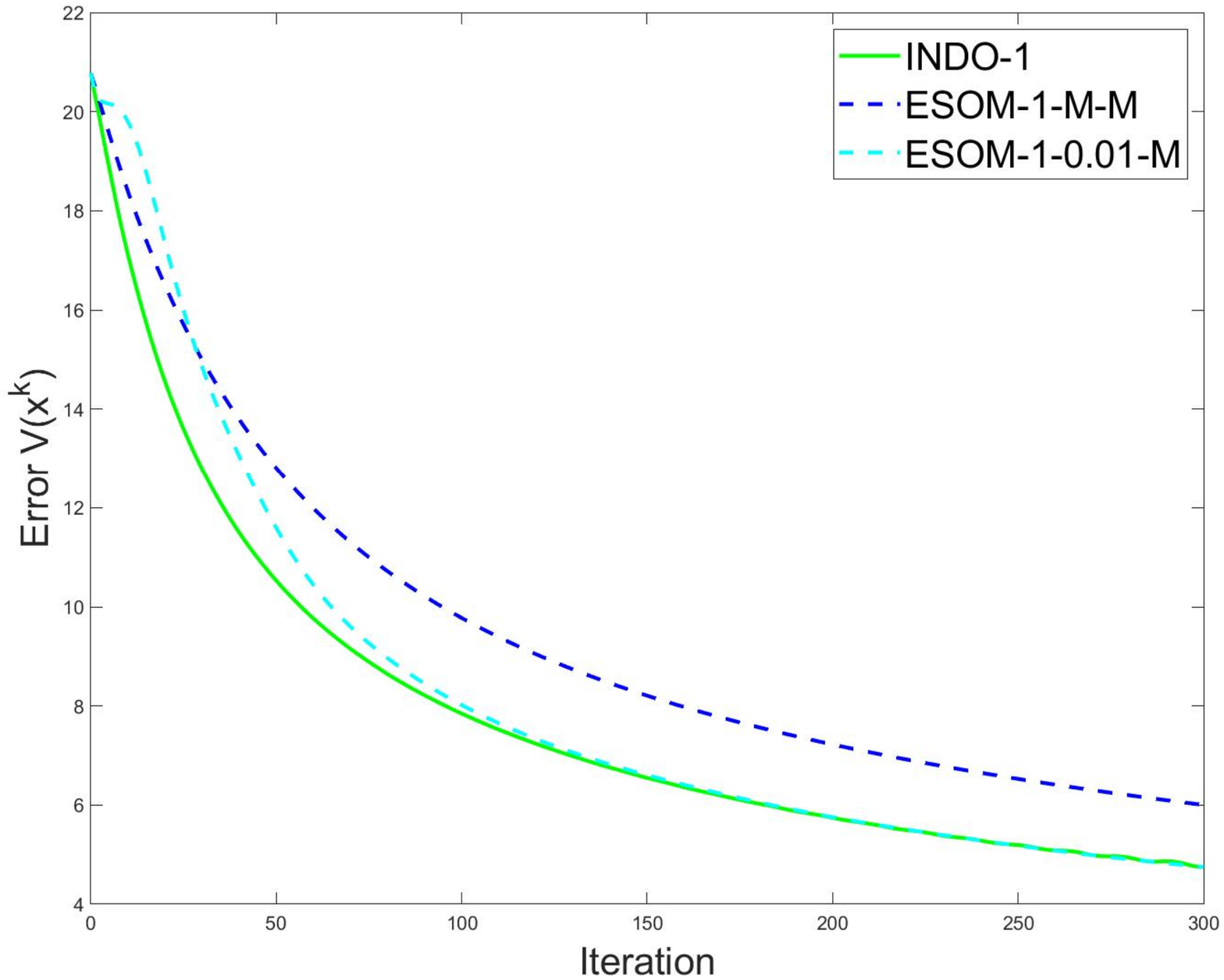}
\includegraphics[width=5.5cm, height=9cm, angle=-90]{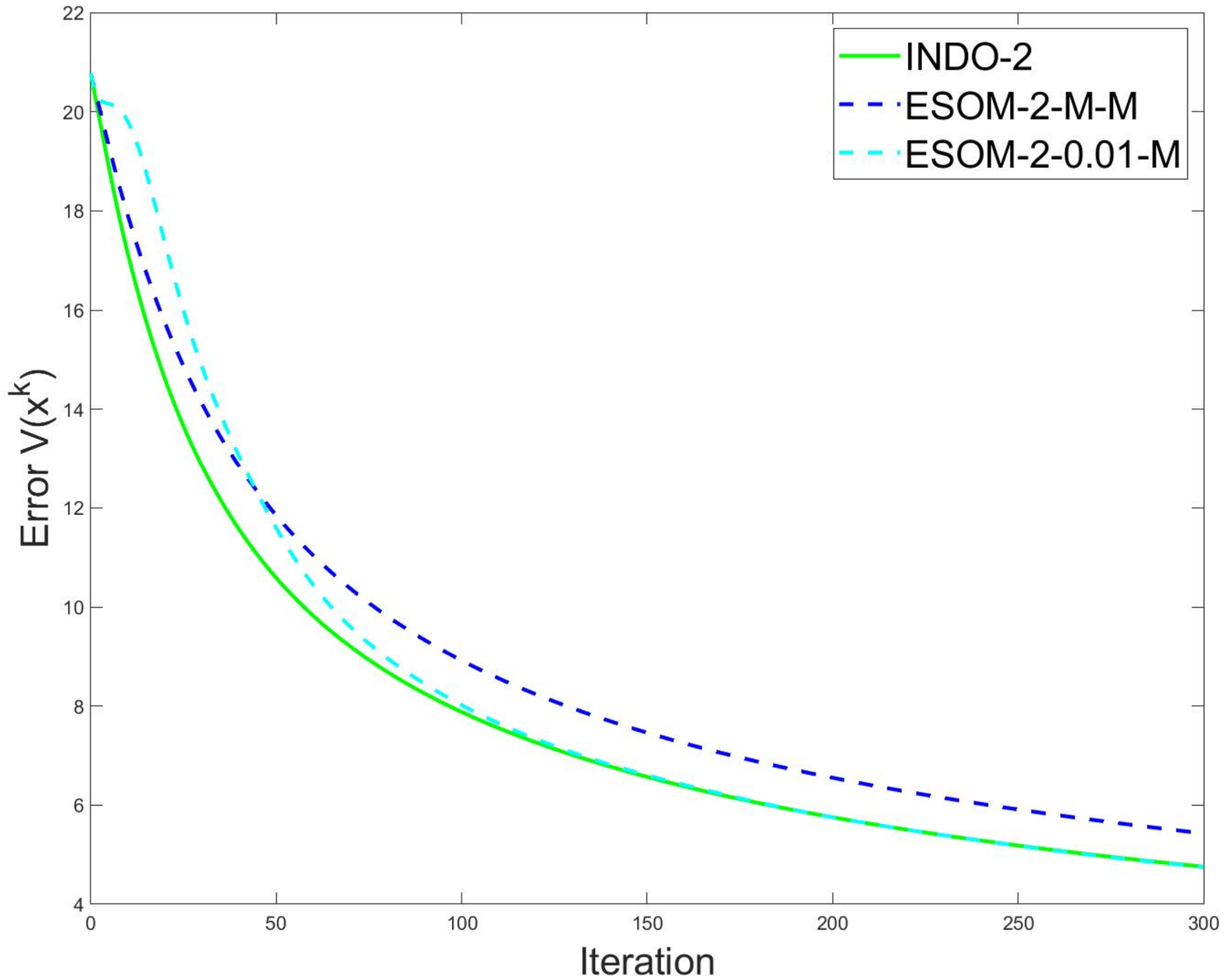}
\includegraphics[width=5.5cm, height=9cm, angle=-90]{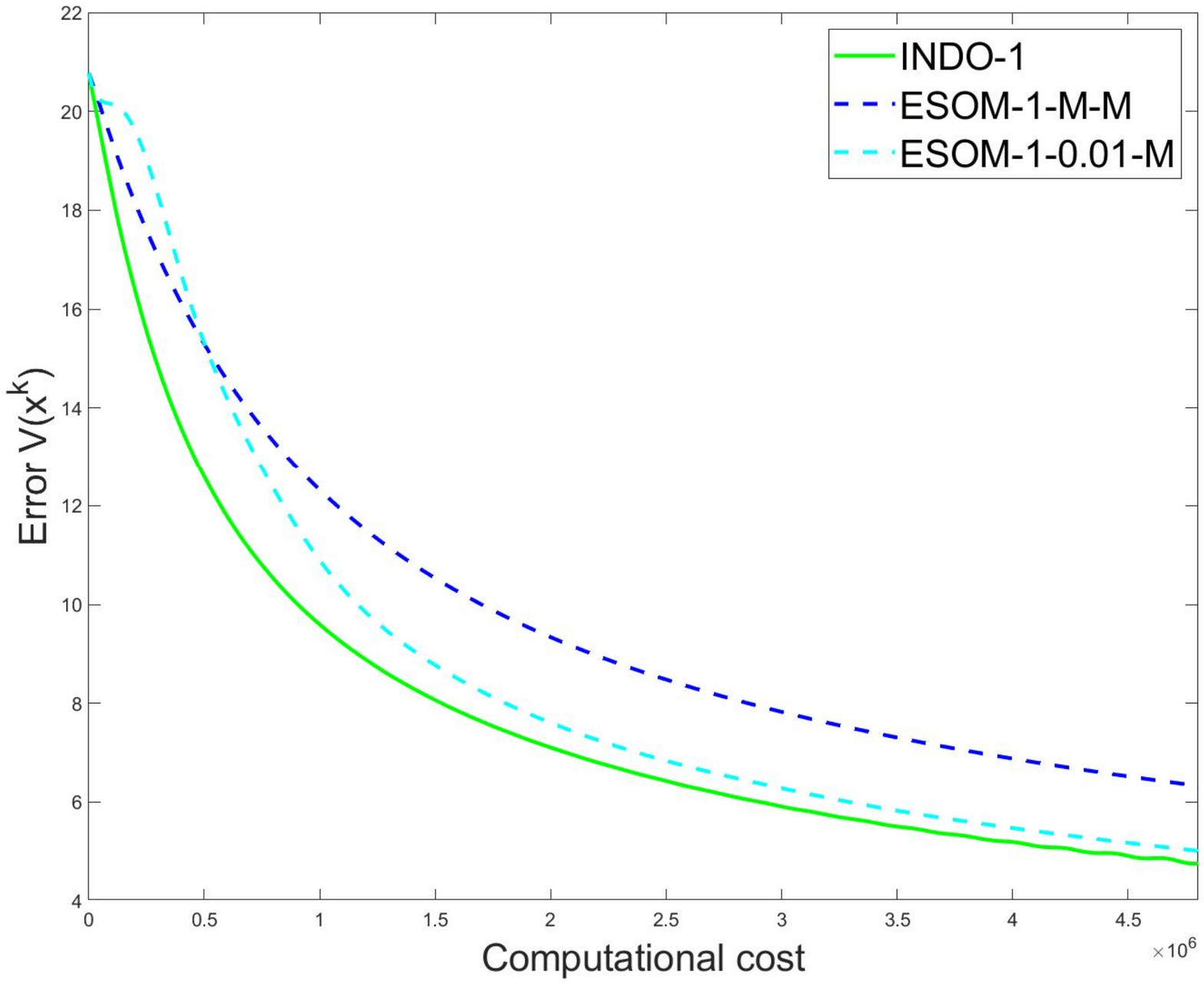}
\includegraphics[width=5.5cm, height=9cm, angle=-90]{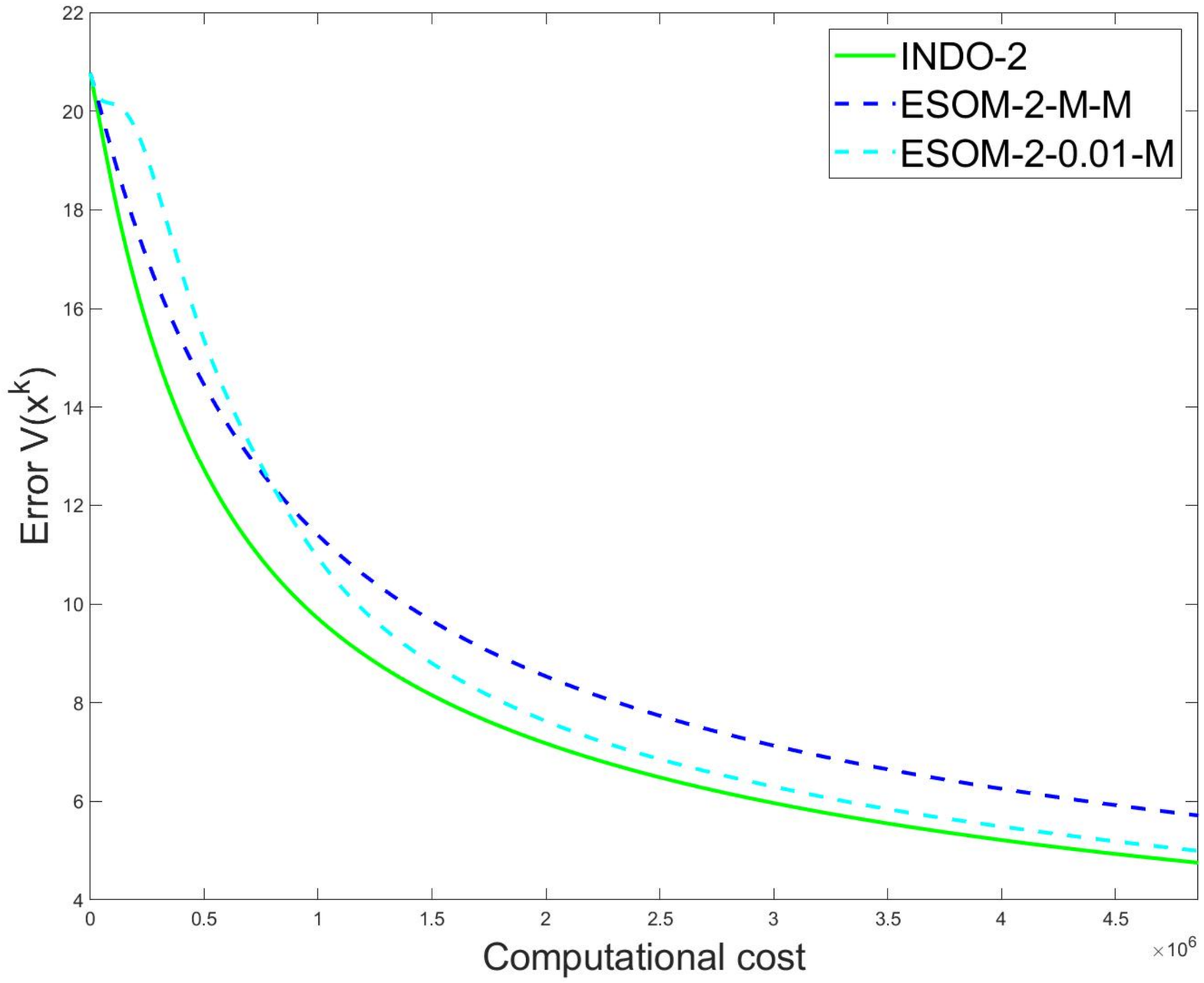}
  \caption{INDO (solid line) versus ESOM methods (dotted lines): error \eqref{errorv} with respect to iterations/communication cost (the two top figures) and computational cost (the two bottom figures) for the number of inner iterations  $\ell=1$ (first and third figure from top) and $\ell=2$ (second and fourth figure from top).  Mushrooms dataset. }
\end{figure}

\begin{figure}
   \includegraphics[width=5.5cm, height=9cm, angle=-90]{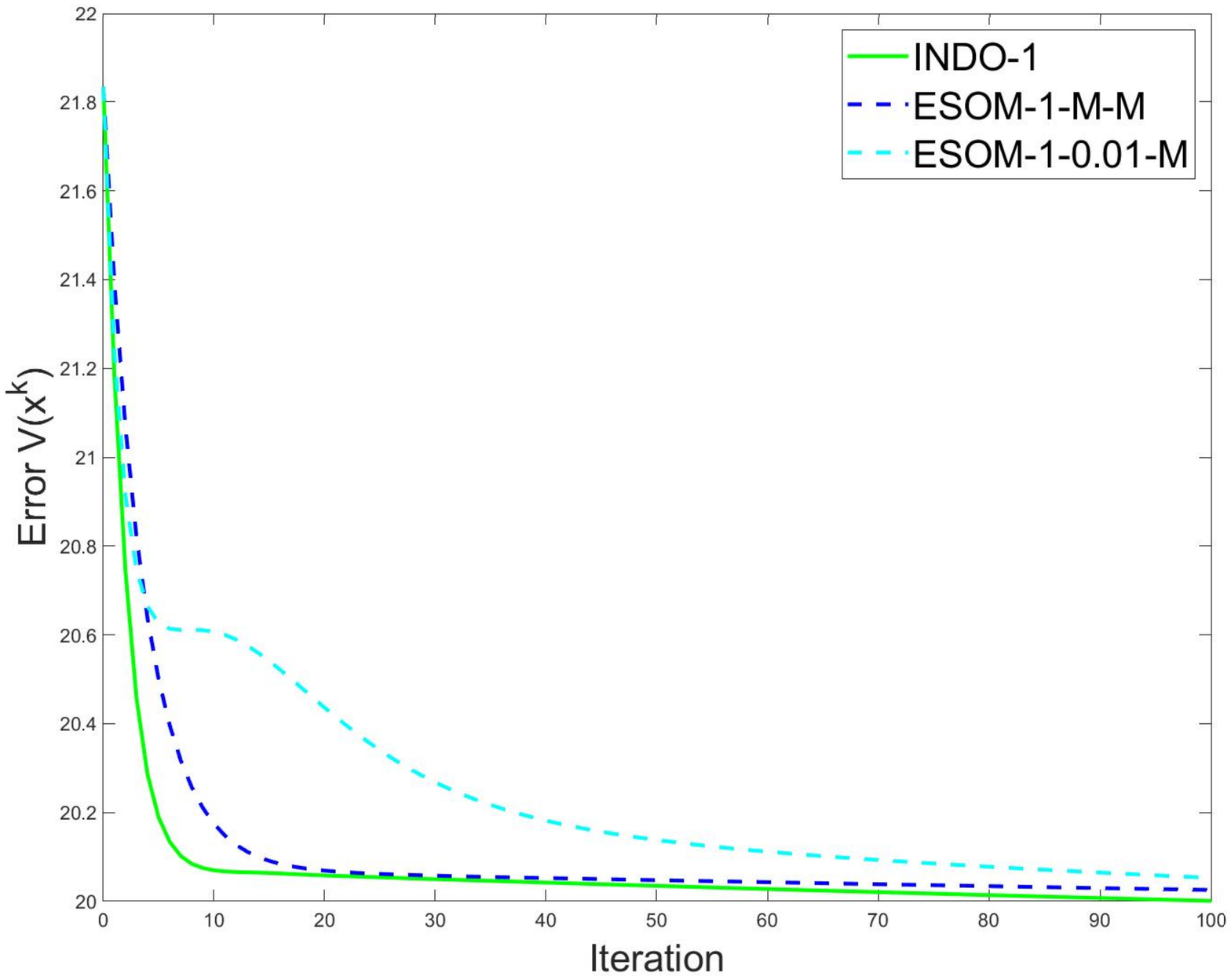}
   \includegraphics[width=5.5cm, height=9cm, angle=-90]{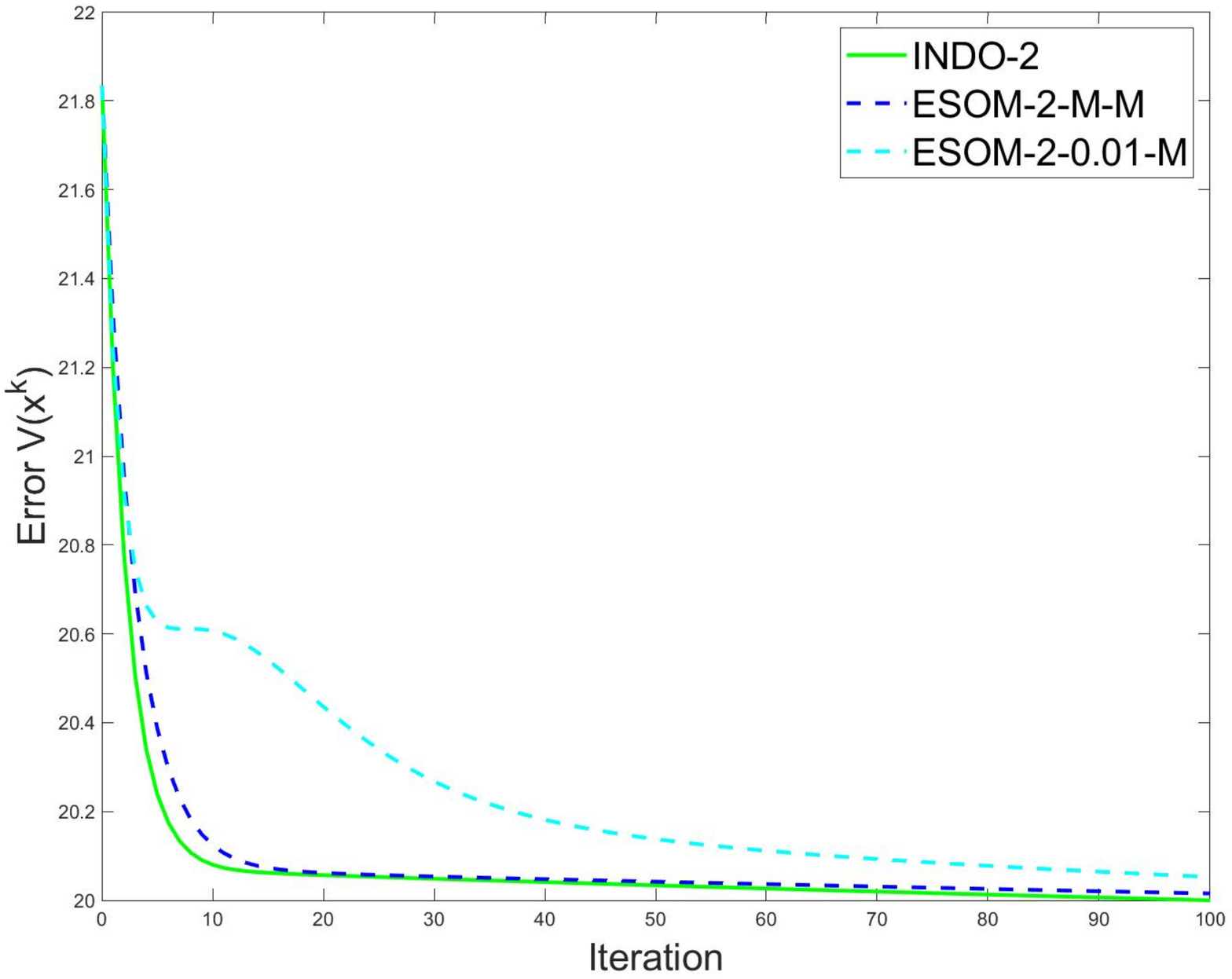}
    \includegraphics[width=5.5cm, height=9cm, angle=-90]{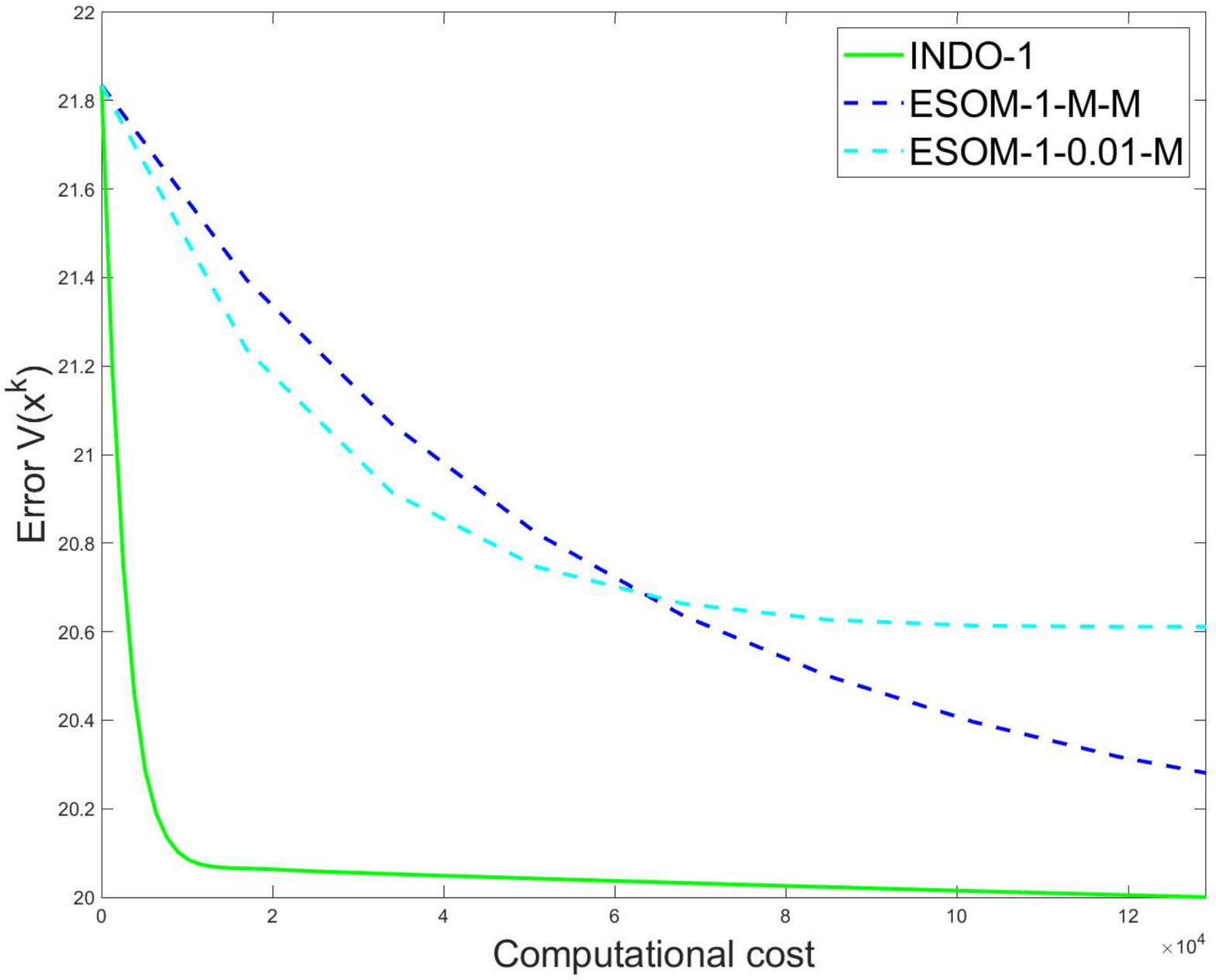}
    \includegraphics[width=5.5cm, height=9cm, angle=-90]{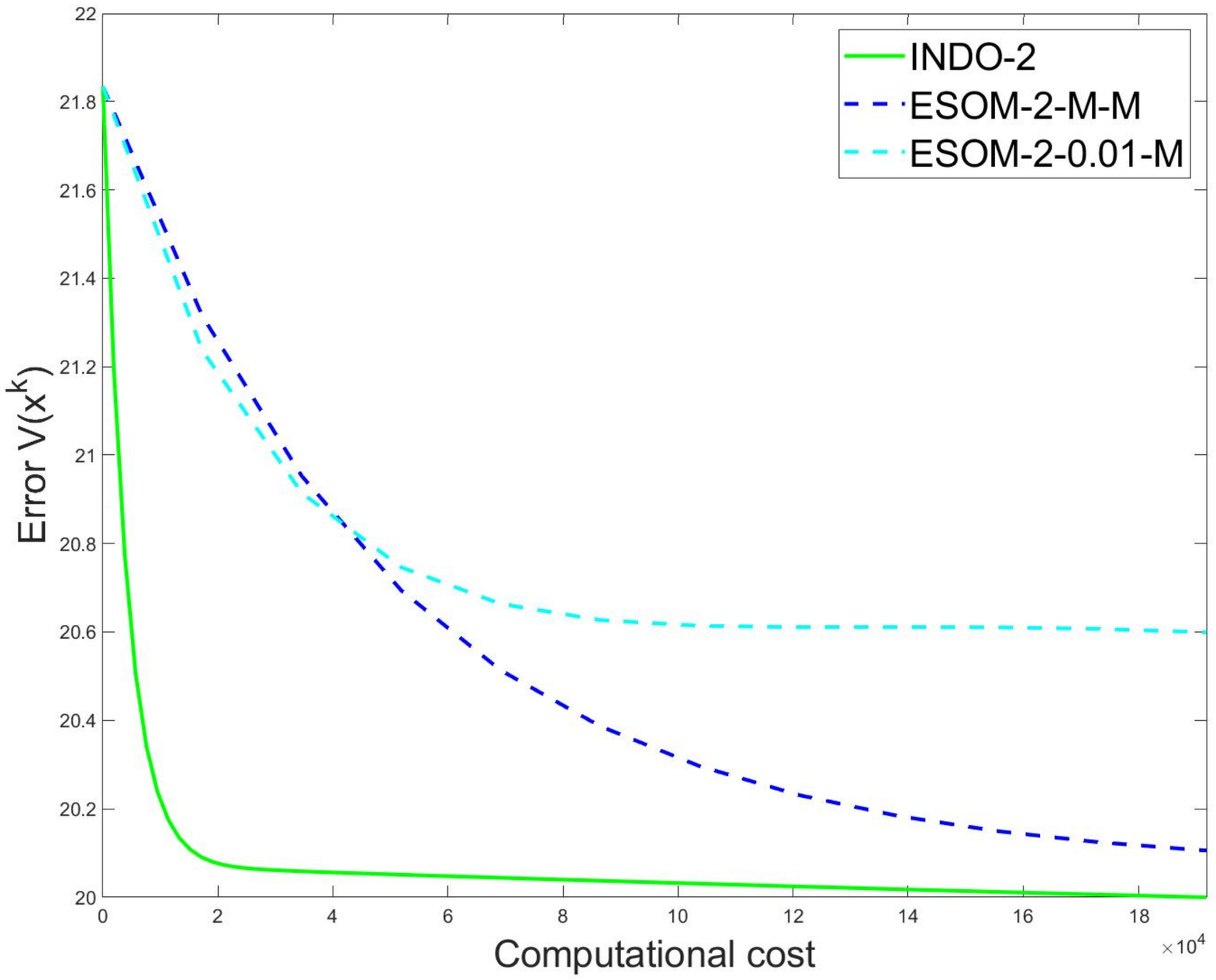}
    \caption{{\footnotesize{INDO (solid line) versus ESOM methods (dotted lines): error \eqref{errorv} with respect to iterations/communication cost (the two top figures) and computational cost (the two bottom figures) for the number of inner iterations  $\ell=1$ (first and third figure from top) and $\ell=2$ (second and fourth figure from top).   LSVT Voice Rehabilitation dataset. }} }
\end{figure}

\begin{figure}
    \includegraphics[width=5.5cm, height=9cm, angle=-90]{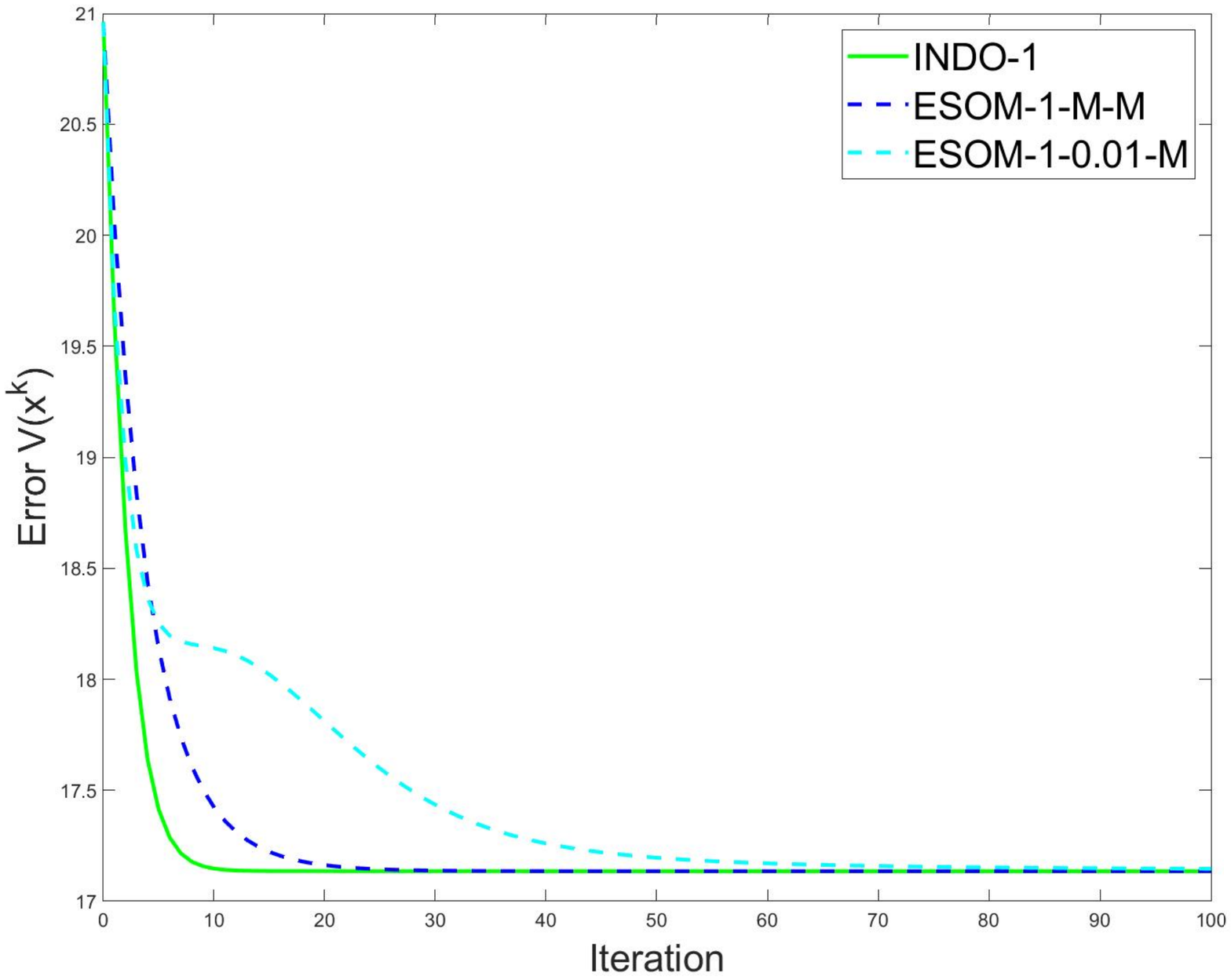}
    \includegraphics[width=5.5cm, height=9cm, angle=-90]{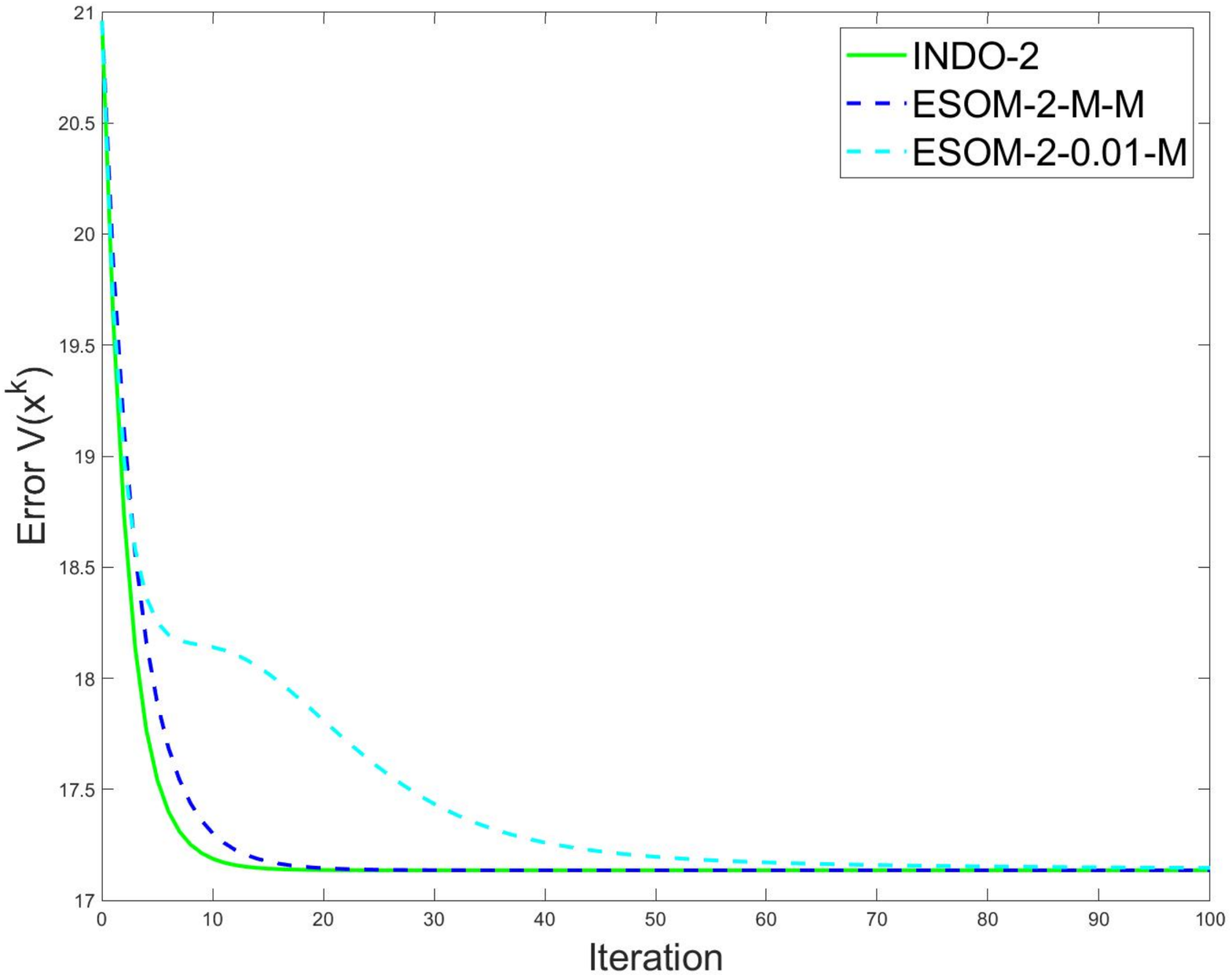}
    \includegraphics[width=5.5cm, height=9cm, angle=-90]{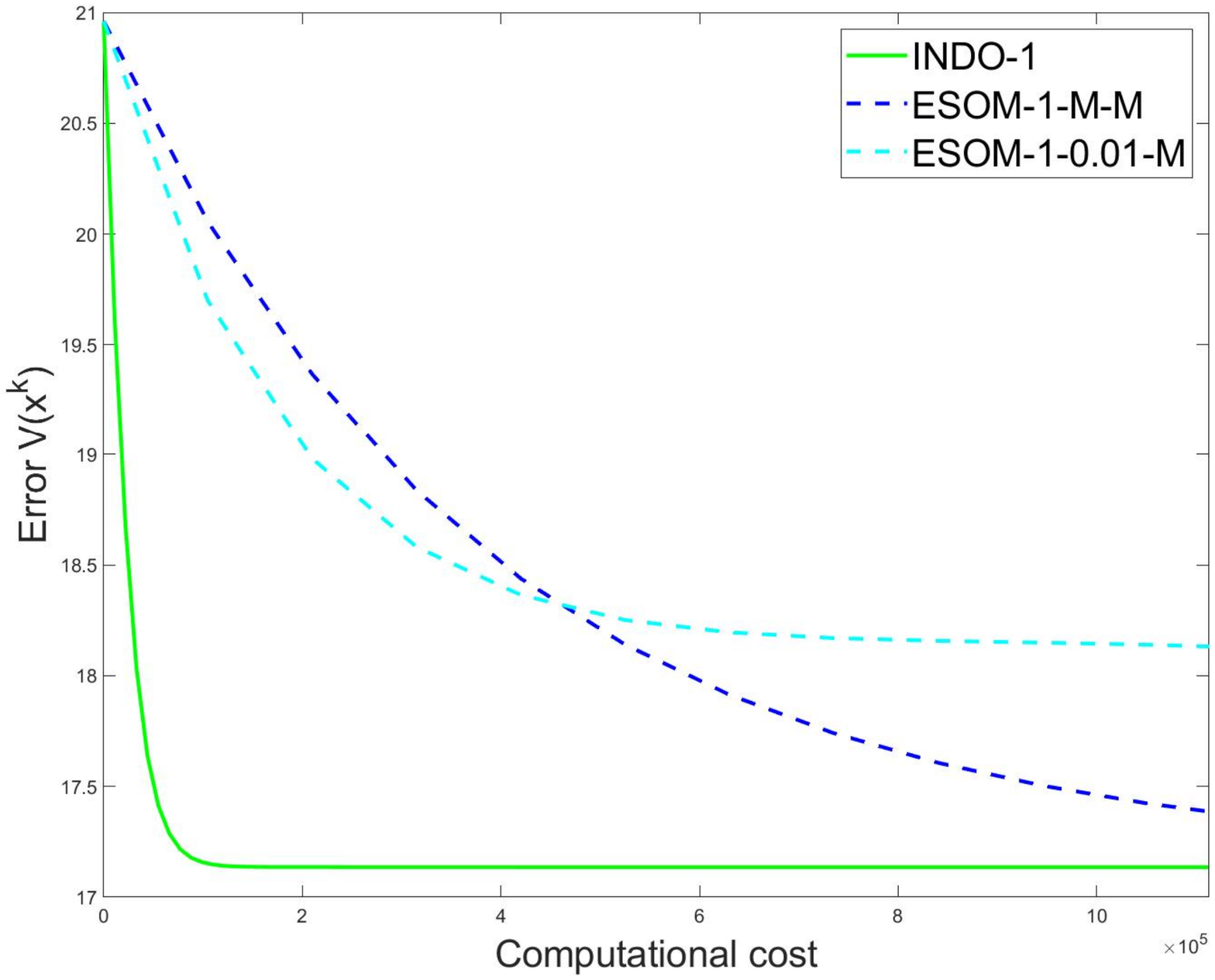}
    \includegraphics[width=5.5cm, height=9cm, angle=-90]{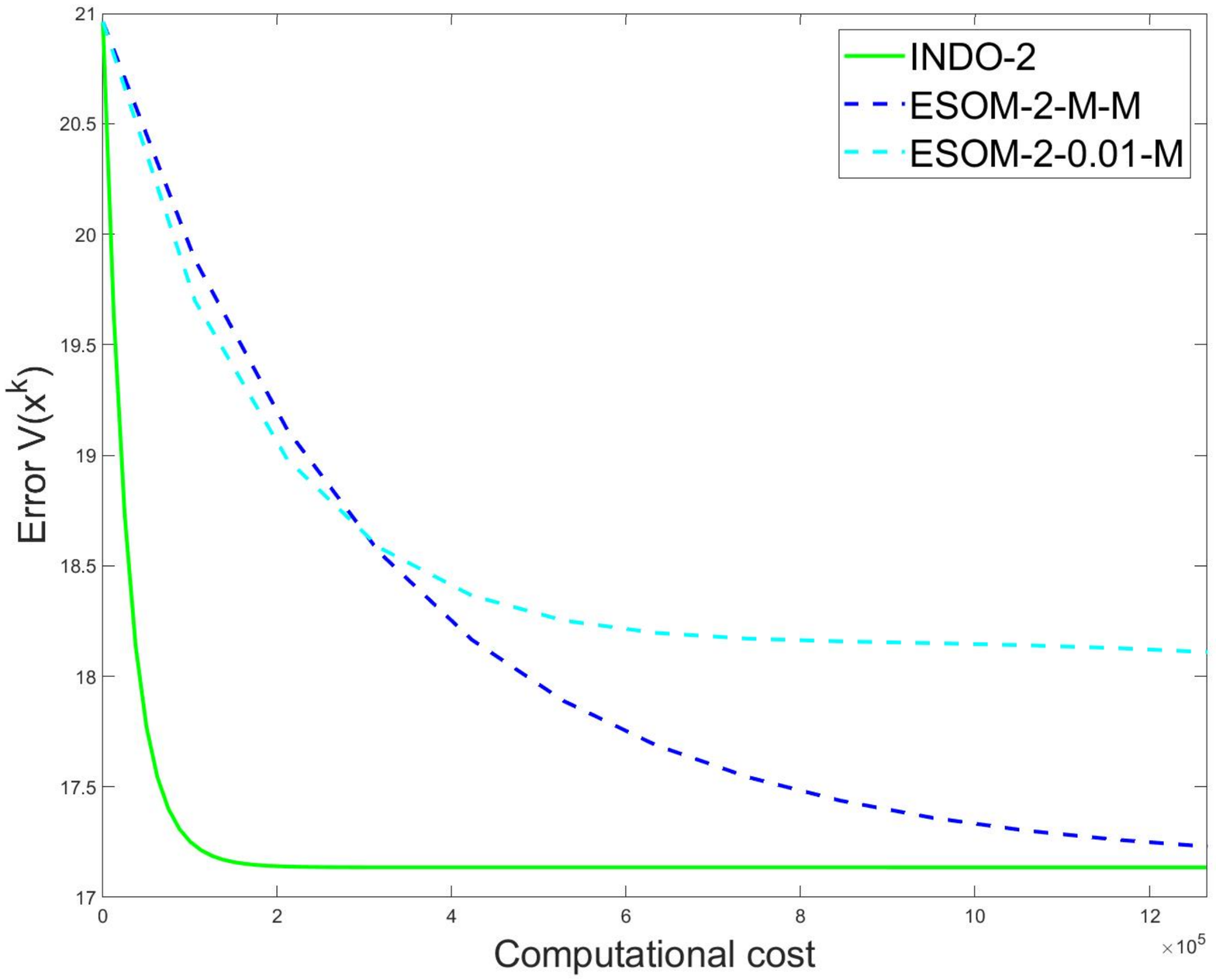}
    \caption{{\footnotesize{INDO (solid line) versus ESOM methods (dotted lines): error \eqref{errorv} with respect to iterations/communication cost (the two top figures) and computational cost (the two bottom figures) for the number of inner iterations  $\ell=1$ (first and third figure from top) and $\ell=2$ (second and fourth figure from top).   Parkinson's Disease Classification dataset.  }} }
\end{figure}

Figures 2-4 present analogous tests on logistic regression problems and the data sets Mushrooms,  LSVT Voice Rehabilitation and Parkinson's Disease Classification. As already mentioned, for these problems $M=1.0001$ so we omit the case $\alpha=0.01, \varepsilon=1$ for ESOM method in these tests. All the parameters for INDO algorithm are the same as in quadratic case.  The results show that INDO algorithm is better or comparable with the ESOM algorithm with respect to iterations, i.e., communication cost. The difference seems to be bigger in the case of $\ell=1$ than for the larger number of inner iterations $\ell=2$. On the other hand, INDO outperforms ESOM method with respect to computational cost.  Notice that for the  Mushrooms dataset,  $n$ is relatively small with respect to $|J_i|$ and thus the common cost is the dominant cost, so the computational cost plots are very similar (in the sense of relative comparison INDO versus ESOM) to the ones that represent iterations and communication costs. For the remaining logistic regression problems, $n$ is comparable to the total sample size $T$ and thus the INDO's advantage in computational cost savings is more evident.

\section{Conclusions}

In this paper, we proposed INDO, an exact distributed second order method for strongly convex distributed consensus optimization. INDO is design based on the framework of Proximal Method of Multipliers (PMM) and maintains a primal and a dual variable per iteration. Unlike existing exact distributed second order methods like ESOM \cite{ESOM} and DQM \cite{DQM}, INDO does not involve explicit Hessian inverse calculation when calculating the primal variable update. Instead, INDO calculates the primal variable update via a few inner iterations of a fixed point method (e.g., Jacobi overrelaxation) to approximately solve the system of linear equations that underlie the Newton direction evaluation; this update process alleviates the need for Newton inverse calculation. We prove that INDO achieves a global linear convergence to the exact solution of the problem of interest. Then, we provide analysis that reveals how the degree of inexactness in solving the Newton direction systems of linear equations affects the overall convergence rate. Furthermore, we provide intuitive bounds on the convergence factor of the involved linear system solver; these bounds shed light on the role of different method’s parameters and provide guidelines on how these parameters should be set. Numerical experiments on several real data sets demonstrates that INDO achieves a comparable speed iteration-wise and communication cost-wise as ESOM, while at the same time reducing computational cost by at least an order of magnitude, when the dimension of the optimization variable is on the order of couple of hundreds or larger.

 \section*{Appendix}
 {\bf Proof of Lemma \ref{L5}}. 
 By assumption A\ref{A2} we have that the aggregate function $ f $ is strongly convex with constant $ m $ and its gradient $ \nabla f $ is Lipschitz continuous with constant $ M. $ Therefore the following inequality holds
\begin{eqnarray}  \label{71}
&\, & \frac{mM}{m+M} \|{\mx}^{k+1}-{\mx}^*\|^2 + \frac{1}{m+M}\|\nabla f({\mx}^{k+1})\\
&-&\nabla f({\mx}^*)\|^2 \nonumber \\
& \leq & ({\mx}^{k+1}-{\mx}^*)^T(\nabla f({\mx}^{k+1}) -\nabla f({\mx}^*)).  \nonumber 
\end{eqnarray}
On the other hand from Lemma \ref{L2} we have 
$$  \nabla f({\mx}^{k+1})- \nabla f({\mx}^*) = -\varepsilon \md^k -  (\II-\WW)^{1/2}({\mv}^{k+1}-{\mv}^*)- \me^k $$
Putting the right hand side of this inequality back in (\ref{71}) and multiplying with $ 2 \alpha $ yields
\begin{eqnarray}  \label{72}
& & \frac{2 \alpha m M}{m+M} \|{\mx}^{k+1}-{\mx}^*\|^2 \\
&+& \nonumber  \frac{2 \alpha}{m+M}\|\nabla f({\mx}^{k+1})-\nabla f({\mx}^*)\|^2 \\
& \leq & - 2 \alpha ({\mx}^{k+1}-{\mx}^*)^T  (\II-\WW)^{1/2}({\mv}^{k+1}-{\mv}^*) \nonumber \\
&-  & 2 \alpha \varepsilon ({\mx}^{k+1}-{\mx}^*)^T \md^k - 2 \alpha ({\mx}^{k+1}-{\mx}^*)^T \me^k.  \nonumber
\end{eqnarray}
Now, Lemma \ref{L1} implies $ \alpha ({\mx}^{k+1}-{\mx}^*)^T  (\II-\WW)^{1/2} = ({\mv}^{k+1}-{\mv}^k)^{T}, $ and therefore
\begin{eqnarray}  \label{73}
& & \frac{2 \alpha m M}{m+M} \|{\mx}^{k+1}-{\mx}^*\|^2 \\
\nonumber 
&+& \frac{2 \alpha}{m+M}\|\nabla f({\mx}^{k+1})-\nabla f({\mx}^*)\|^2 \\
& \leq & - 2({\mv}^{k+1}-{\mv}^k)^T({\mv}^{k+1}-{\mv}^*) - 2 \alpha \varepsilon ({\mx}^{k+1}-{\mx}^*)^T \md^k \nonumber \\
&-& 2 \alpha ({\mx}^{k+1}-{\mx}^*)^T \me^k.  \nonumber
\end{eqnarray}
Notice further that
\begin{eqnarray*}
 &\,&2({\mv}^{k+1}-{\mv}^k)^T({\mv}^{k+1}-{\mv}^*) = \|{\mv}^{k+1} - {\mv}^k\|^2 \\
 &+& \|{\mv}^{k+1}-{\mv}^*\|^2 - \|{\mv}^k - {\mv}^*\|^2 
 \end{eqnarray*}
and
$$ 2 ({\mx}^{k+1}-{\mx}^*)^T \md^k = \|\md^k\|^2 + \|{\mx}^{k+1}-{\mx}^*\|^2 - \|{\mx}^k-{\mx}^*\|^2. $$ 
Going back to (\ref{73}) we get
\begin{eqnarray}  
\nonumber 
& & \frac{2 \alpha m M}{m+M} \|{\mx}^{k+1}-{\mx}^*\|^2 + \frac{2 \alpha}{m+M}\|\nabla f({\mx}^{k+1})-\nabla f({\mx}^*)\|^2 \\
& \leq & - \|{\mv}^{k+1} - {\mv}^k\|^2 - \|{\mv}^{k+1}-{\mv}^*\|^2 + \|{\mv}^k - {\mv}^*\|^2 \label{74} \\
&-& \alpha \varepsilon \|\md^k\|^2 - \alpha \varepsilon \|{\mx}^{k+1}-{\mx}^*\|^2 \nonumber  \\
& + & \alpha \varepsilon \|{\mx}^k-{\mx}^*\|^2 - 2 \alpha ({\mx}^{k+1}-{\mx}^*)^T \me^k.  \nonumber
\end{eqnarray}
Lemma \ref{L1} implies that $ \|{\mv}^{k+1} - {\mv}^k\|^2 = \|{\mx}^{k+1}-{\mx}^*\|^2_{\alpha^2(\II-\WW)}$ and by definition of Lyapunov function we have 
$ \|\uu^k-\uu^*\|^2_{\cal G} - \|\uu^{k+1} - \uu^*\|^2_{\cal G} =  \|{\mv}^k - {\mv}^*\|^2 - \|{\mv}^{k+1} - {\mv}^*\|^2+\alpha \varepsilon \|{\mx}^k-{\mx}^*\|^2   - \alpha \varepsilon \|{\mx}^{k+1}-{\mx}^*\|^2. $ Thus, inequality (\ref{75}) reduces to 
\begin{eqnarray}  \label{74}
& & \frac{2 \alpha m M}{m+M} \|{\mx}^{k+1}-{\mx}^*\|^2 \\
&+& \nonumber  \frac{2 \alpha}{m+M}\|\nabla f({\mx}^{k+1})-\nabla f({\mx}^*)\|^2 \\
& \leq &  \|\uu^k-\uu^*\|^2_{\cal G} - \|\uu^{k+1} - \uu^*\|^2_{\cal G}  - \alpha \varepsilon \|\md^k\|^2 \nonumber \\
&-&  \|{\mx}^{k+1}-{\mx}^*\|^2_{\alpha^2(\II-\WW)} \nonumber \\ &- & 2 \alpha ({\mx}^{k+1}-{\mx}^*)^T \me^k. \nonumber
\end{eqnarray}
Regrouping the terms in the above inequality yields
\begin{eqnarray}  \label{75}
& & \|\uu^{k+1}-\uu^*\|^2_{\cal G} - \|\uu^{k} - \uu^*\|^2_{\cal G} \\
& \leq &- \frac{2 \alpha}{m+M}\|\nabla f({\mx}^{k+1}) \nonumber \\
&-&\nabla f({\mx}^*)\|^2  - \|{\mx}^{k+1}-{\mx}^*\|^2_{\frac{2 \alpha m M}{m+M} \II+\alpha^2(\II-\WW)} \nonumber \\
& -&  \alpha \varepsilon \|\md^k\|^2 - 2 \alpha ({\mx}^{k+1}-{\mx}^*)^T \me^k. \nonumber 
\end{eqnarray}
Using the bound $ 2 ({\mx}^{k+1}-{\mx}^*)^T \me^k \geq -1/\zeta\|{\mx}^{k+1}-{\mx}^*\|^2-\zeta \|\me^k\|^2, $ which holds for any $ \zeta> 0, $ we get
 \begin{eqnarray} 
 & & \|\uu^{k+1}-\uu^*\|^2_{\cal G} - \|\uu^{k} - \uu^*\|^2_{\cal G} \\
& \leq &- \frac{2 \alpha}{m+M}\|\nabla f({\mx}^{k+1})-\nabla f({\mx}^*)\|^2  - \|{\mx}^{k+1} \nonumber \\
&-&{\mx}^*\|^2_{(\frac{2 \alpha m M}{m+M} -\frac{\alpha}{\zeta})\II+\alpha^2(\II-\WW)} \nonumber \\
& -&  \alpha \varepsilon \|\md^k\|^2 +  \alpha \zeta \|\me^k\|^2. \nonumber
\end{eqnarray}

\begin{thebibliography}{99}
 
\bibitem{koreni}  Bai, Z., Silverstein,  J. W.,  Spectral Analysis of Large Dimensional Random Matrices, Springer, 2010.
 
 \bibitem{DQN} Bajovi\'c, D., Jakoveti\'c,D.,  Kreji\'c, N., Krklec Jerinki\'c, N., Newton-like Method with Diagonal Correction for Distributed Optimization,  SIAM Journal on Optimization, 27,  2 (2017), 1171-1203
 
 

\bibitem{novo3}  Baingana, B.,  Giannakis, G., B.,  Joint Community and Anomaly Tracking in Dynamic Networks,  IEEE Transactions on Signal Processing, 64(8), (2016), pp. 2013-2025.

\bibitem{IMA} Bellavia, S., Kreji\'c, N., Krklec Jerinki\'c, N., Subsampled Inexact Newton Methods for minimizing large sums of convex functions, IMA J. Numer. Anal. 40,4 (2020), 2309-2341.

\bibitem{wei} Berahas, A. S.,  Bollapragada, R.,   Keskar, N. S.,   Wei, E., Balancing Communication and Computation in Distributed Optimization, IEEE Transactions on Automatic Control, 64(8), (2019), pp. 3141-3155.



 \bibitem{BoydADMM}
Boyd, S., Parikh, N., Chu, E., Peleato, B., Eckstein, J., Distributed   optimization and statistical learning via the alternating direction method of   multipliers, Foundations and Trends in Machine Learning, 3(1),  (2011) pp. 1-122.

\bibitem{SayedEstimation}
Cattivelli, F.,  Sayed,  A.~H., Diffusion {LMS} strategies for distributed   estimation, 
 IEEE Transactions on  Signal  Processing,  58(3),  (2010) pp.   1035--1048.

\bibitem{EISEN} Eisen, M., Mokhtari, A., Ribeiro, A., 
Decentralized Quasi-Newton Methods, IEEE Trans. Signal Process, vol. 65, no. 10, pp. 2613--2628, 2017.

\bibitem{DES} Dembo, R.S., , Eisenstadt, S.C., Steihaugh, T., Inexact Newton Methods, SIAM Journal on Numerical Analysis, 19,2, (1982), 400-408. 



\bibitem{r1} Eisen, M., Mokhtari, A., Ribeiro, A., A Primal-Dual Quasi-Newton Method for Exact Consensus Optimization, IEEE Transactions on Signal Processing, Vol. 67, No. 23, Dec. 2019., pp. 5983--5997.
%
\bibitem{AG} Greenbaum, A.,  Iterative Methods for Solving Linear Systems, SIAM, 1997.

\bibitem{dusan} Jakoveti\'c, D., A Unification and Generalization of Exact Distributed First Order Methods, IEEE Transactions on Signal and Information Processing over Networks, 5(1), (2019), pp. 31-46.

\bibitem{dragana} Jakoveti\'c, D., Bajovi\'c, D., Xavier, J., Moura, J.M.F., Primal-Dual Methods for Large-Scale and Distributed Convex Optimization and Data Analytics, Procedings of the IEEE, 108,11 (2020), 1923-1938. 
 
 
 \bibitem{DFIX} Jakoveti\'c, D.,  Kreji\'c, N.,  Krklec Jerinki\'c, N. , Malaspina, G. ,  Micheletti, A.,  Distributed Fixed Point Method for Solving Systems of Linear Algebraic Equations,
	Automatica (2021), Vol. 134.

 


  \bibitem{arxivVersion}
Jakoveti\'c, D., Xavier, J.,  Moura, J.~M.~F.,  Fast distributed gradient
  methods, \emph{IEEE Transactions on  Automatic Control}, 59(5), (2014)  pp. 1131--1146.

\bibitem{Espectral} Jakoveti\'c, D.,  Kreji\'c, N.,  Krklec Jerinki\'c, N.,  Exact spectral-like gradient method for distributed optimization, Computational Optimization and Applications, 74,  (2019), pp. 703–728.

\bibitem{EFIX} Jakoveti\'c, D.,  Kreji\'c, N.,  Krklec Jerinki\'c, N.,  EFIX: Exact Fixed Point Methods for Distributed Optimization, arxiv preprint,  arXiv:2012.05466v1  (2020).


  
 \bibitem{novo1}  Lee, J. M.,  Song, I.,  Jung, S.,   Lee, J.,  A rate adaptive convolutional coding method for multicarrier DS/CDMA systems, MILCOM 2000 Proceedings 21st Century Military Communications. Architectures and Technologies for Information Superiority (Cat. No.00CH37155), Los Angeles, CA, (2000), pp. 932-936. 
 
 \bibitem{r2} Mansoori, F., Wei, E., A Fast Distributed Asynchronous Newton-Based Optimization Algorithm,  IEEE Transactions on Automatic Control, Vol. 65, No. 7, July 2020, pp. 2769--2784.
 

\bibitem{NN} Mokhtari, A., Ling, Q., Ribeiro, A., Network Newton Distributed Optimization Methods, IEEE Transactions on Signal Processing, 65,1 (2017), 146-161. 
%
%
\bibitem{ESOM} Mokhtari, A.,  Shi,  W.,   Ling, Q.,    Ribeiro, A.,
 A Decentralized Second Order Method with Exact Linear Convergence Rate for Consensus Optimization,
IEEE Transactions on Signal and Information Processing over Networks, 2(4), (2016), pp. 507-522.

\bibitem{DQM} Mokhtari, A.,  Shi,  W.,   Ling, Q.,    Ribeiro, A.,
 DQM: Decentralized quadratically approximated alternating direction method of multipliers, IEEE Transactions on Signal Processing, vol. 64, no. 19, pp. 5158-5173, Oct. 2016.
%


 \bibitem{JoaoMotaMPC}
Mota, J., Xavier, J., Aguiar, P., P\"uschel, M., Distributed optimization
 with local domains: Applications in MPC and network flows, IEEE Transactions on Automatic Control, 60(7), (2015), pp. 2004-2009. 

%

 
\bibitem{harnessing} Qu, G., Li, N., Harnessing smoothness to accelerate distributed optimization, IEEE Transactions on Control of Network Systems, 5(3), (2018),  pp. 1245-1260. 

\bibitem{PARKINSON} Sakar, C., Serbes, Gorkem, Gunduz, Aysegul, Nizam, Hatice, Sakar, Betul. (2018). Parkinson's Disease Classification. UCI Machine Learning Repository.

%
\bibitem{scutari2} Scutari, G., Sun, Y., Parallel and Distributed Successive Convex Approximation Methods for Big-Data Optimization,
 Multi-agent Optimization, Lecture Notes in Mathematics, pp. 141-308,  Springer, 2018.


\bibitem{extra}
Shi, W., Ling, Q., Wu, G., Yin, W.,  EXTRA: an Exact First-Order Algorithm for Decentralized Consensus Optimization,
SIAM Journal on Optimization, 2(25), (2015), pp. 944-966.


 
 

\bibitem{d3} Sun, Y., Daneshmand, A., Scutari, G., Convergence Rate of Distributed
Optimization Algorithms based on Gradient Tracking, 
arXiv:1905.02637, (2019).


\bibitem{VOICE} Tsanas, Athanasios. (2014). LSVT Voice Rehabilitation. UCI Machine Learning Repository


\bibitem{d5} Tian, Y., Sun, Y., Scutari, G., Achieving Linear Convergence in Distributed Asynchronous Multi-agent Optimization, IEEE Trans. on
Automatic Control, (2020).

\bibitem{d6} Tian, Y., Sun, Y., Scutari, G., Asynchronous Decentralized Successive
Convex Approximation, arXiv:1909.10144, (2020).
%

\bibitem{mush} UCI Machine Learning Expository, https://archive.ics.uci.edu/ml/datasets/Mushroom.


\bibitem{d8} Xin, R.,  Khan, U. A., Distributed Heavy-Ball: A Generalization and Acceleration of First-Order Methods With Gradient Tracking,  IEEE Transactions on Automatic Control, 65(6), (2020), pp. 2627-2633.


 \bibitem{optra} Xu, J., Tian, Y., Sun, Y., Scutari G., Accelerated primal-dual algorithms for distributed smooth convex optimization over networks,  International Conference on Artificial Intelligence and Statistics, PMLR, (2020), pp. 2381-2391.


\bibitem{d10} Xu, J., Tian, Y., Sun, Y., Scutari G., Distributed Algorithms for Composite Optimization: Unified Framework and Convergence Analysis,  IEEE Transactions on Signal Processing, Vol. 69, pp. , June 2021.3555 - 3570, June 2021.

\bibitem{Nedickaskadni}  Yousefian, F.,  Nedi\'c, A.,  Shanbhag, U. V., On stochastic gradient and subgradient methods with adaptive steplength sequences, Automatica, 48(1), (2012), pp. 56-67.


 \bibitem{ypma} Ypma, T., Local convergence of Inexact Newton Methods, SIAM Journal on Numerical Analysis, 21(3), (1984), 583–590. 
 
 \bibitem{Sayed2} Yuan, K., Ying, B., Zhao, X., Sayed, A. H., Exact diffusion for
distributed optimization and learning -- Part I: Algorithm development,  IEEE Transactions on Signal Processing, Vol. 67, No. 3, Feb., 2019., pp. 708--723.

\bibitem{Sayed3} Yuan, K., Ying, B., Zhao, X., Ali H. Sayed, 
Exact Diffusion for Distributed Optimization and Learning -- Part II: Convergence Analysis, IEEE Trans. Signal Process., vol. 67, No. 3, pp. 724--739, 2019.

\bibitem{UsmanXin}
Xin, R., Khan, U. A., and Kar, S., Variance-reduced decentralized stochastic optimization with accelerated convergence, IEEE Transactions on Signal Processing, vol. 68, pp. 6255--6271, Oct. 2020.

 \bibitem{Zhang} Zhang, J., Ling, Q., So, A.-M., A Newton Tracking Algorithm with Exact Linear Convergence Rate for Decentralized Consensus Optimization, IEEE Transactions on Signal and Information Processing over Networks, 7, (2021), pp. 346--358. 
 
  \bibitem{Max}
 Zhang, S., Tepedelenlioglu, C., Banavar, M.K., Spanias, A., Max Consensus in Sensor Networks: Non-Linear Bounded Transmission and Additive Noise,
  IEEE Sensors Journal, Vol. 16, No. 24, pp.  9089--9098, Dec. 2016.








  
%
%
%
%
%

%
%

%



%
%
%




%
%


%
%










%

%
%
%



 \end{thebibliography}
\end{document}